\documentclass[reqno,12pt]{amsart} 
\usepackage{amsmath}
\usepackage{amssymb}
\usepackage{amsthm}
\newcommand\NoBlackBoxes{\global\overfullrule0pt}
\NoBlackBoxes
\theoremstyle{plain} 
\newtheorem{theorem}{Theorem} 

\newtheorem{lemma}[theorem]{Lemma}
\newtheorem{proposition}[theorem]{Proposition}

\newtheorem{corollary}[theorem]{Corollary}
\def\4{\kern1pt}

\def\6{\vphantom0}

\def\8{\kern-10pt}
\def\7#1{_{(#1)}}

\theoremstyle{definition}
\newtheorem*{assumption (H1)}{Assumption (H1)}
\newtheorem*{assumption (H2)}{Assumption (H2)}

\theoremstyle{remark}
\newtheorem{remark}[theorem]{Remark}

\numberwithin{equation}{section}
\numberwithin{theorem}{section}
\makeatletter
\let\serieslogo@\relax
\let\@setcopyright\relax

\def\speciallabelmark#1{\def\@currentlabel{#1}}
\makeatother

\renewcommand{\Bbb}{\mathbb}

\newcommand{\Cal}{\mathcal}

\newcommand{\beqa}{\begin{eqnarray}}
\newcommand{\beqan}{\begin{eqnarray*}}
\newcommand{\eeqa}{\end{eqnarray}}
\newcommand{\eeqan}{\end{eqnarray*}}
\def\beq#1\eeq{\begin{equation}#1\end{equation}}

\begin{document}

\def\ffrac#1#2{\raise.5pt\hbox{\small$\4\displaystyle\frac{\,#1\,}{\,#2\,}\4$}}
\def\ovln#1{\,{\overline{\!#1}}}
\def\ve{\varepsilon}
\def\kar{\beta_r}

\title{Limit theorems in free probability theory. I}

\author{G. P. Chistyakov$^{1}$}
\address
{Gennadii Chistyakov, 
Institute for Low Temperature Physics and Engineering,\newline
National Academy of Sciences of Ukraine,
47 Lenin Ave.,
61103 Kharkov, 
Ukraine \smallskip}
\email {chistyakov@ilt.kharkov.ua, \,chistyak@mathematik.uni-bielefeld.de} 

\author{F. G\"otze$^{1}$}
\thanks{1) Research supported by the DFG-Forschergruppe 399/2 and SFB 701. 
Partially supported by INTAS grant N 03-51-5018.}
\address
{Friedrich G\"otze,
Fakult\"at f\"ur Mathematik,
Universit\"at Bielefeld,
Postfach 100131,\newline
33501 Bielefeld, 
Germany}
\email {goetze@mathematik.uni-bielefeld.de}



\begin{abstract}
Based on a new analytical approach to the definition of additive free
convolution on 
probability measures on the
real line we prove free analogs of limit theorems for sums
for non-identically distributed random variables
in classical Probability Theory.
\end{abstract}

\maketitle
\markboth{ G. P. Chistyakov and F. G\"otze}{Limit Theorems
in Free Probability Theory}

\section{Introduction}
In recent years a number of papers are investigating limit theorems
for the free convolution of probability measures (p-measures)
defined  by D. Voiculescu.
The~key concept of this definition is the~notion of freeness,
which can be interpreted as a~kind of independence for
noncommutative random variables. As in the~classical probability where
the~concept of independence gives rise to the~classical convolution, 
the~concept of freeness leads to a~binary operation on the~p-measures 
on the~real line, the~free convolution. Many classical results 
in the~theory of addition of independent random variables have their
counterpart in this new theory, such as the~law of large numbers,
the~central limit theorem, the~L\'evy-Khintchine formula and others.
We refer to Voiculescu, Dykema and Nica \cite{Vo:1992} for introduction 
to these topics. Bercovici and Pata \cite{BeP:1999} established
the~distributional behavior of sums of free identically distributed 
random variables and described explicitly the~correspondence 
between limits laws for free and classical additive convolution. 
In this paper, using a~new approach
to the~definition of the~additive free convolution (see \cite{ChG:2005}),
we generalize the~results of Bercovici and Pata to the~case of free
non-identically distributed random variables. We show that the~parallelism  
found by Bercovici and Pata holds in the~common case of
free non-identically distributed random variables. 
Our approach to the~definition
of the~additive free convolution allows to obtain estimates of the~rate
of convergence of distribution functions of free sums. We prove
the~semi-circle approximation theorem (an analog of the~Berry-Esseen
inequality), the~law of large number with estimates of the~rate of 
convergence. 
We describe L\'evy's class $\Cal L_{\boxplus}$ 
of limiting distributions of normed sums
of free random variables obeying infinitesimal conditions. 
As in the~classical case we prove the~norming theorem, 
which is necessary and sufficient conditions for convergence,
(see \cite{GKo:1968} and \cite{Lo:1963}) and derive 
the~canonical representation of the~measures of the~class~$\Cal L_{\boxplus}$.
Furthermore, we shall give a~characterization of the~class~$\Cal L_{\boxplus}$
by means of the~property of self-decomposability, extending results 
by Barndorff-Nielsen and Thorbj{\o}rsen~\cite{BarTho:2002}. 

The~paper is organized as follows. In Section~2 we formulate and
discuss the main results of the~paper. In Section~3 we formulate
auxiliary results. In Section~4 we prove the~extended additive free
central limit theorem for general case of free {\it non-identically}
distributed random variables.
This extends the~Bercovici and Pata parallelism
between free additive and classical additive infinite divisibility
and limits laws for free and classical additive convolution to the
general case.
In Section~5, using results of Section~4, we describe an~analog
of the~L\'evy class $\mathcal L_{\boxplus}$ for additive free convolution.
We establish the~Bercovici and Pata parallelism between the~classical
~L\'evy class $\mathcal L$ and the~class $\mathcal L_{\boxplus}$.   
In Section~6, using our approach to the~definition of the~additive
free convolution, we derive the~semicircle approximation theorem
(an~analog of the~Berry-Esseen inequality) as well as a~law of large numbers 
with estimates of convergence.

\section{Results}
Denote by $\mathcal M$ the~family of all Borel p-measures
defined on the~real line $\Bbb R$. On $\mathcal M$ define 
the~associative composition laws denoted $*$ and $\boxplus$ as follows.
For $\mu_1,\mu_2\in\mathcal M$ let the~p-measure $\mu_1*\mu_2$ denote
the~classical convolution of $\mu_1$ and $\mu_2$. 
In probabilistic terms, $\mu_1*\mu_2$
is the~probability distribution of $X+Y$, where $X$ and $Y$ are
(commuting) independent random variables with probability distributions $\mu_1$ 
and $\mu_2$ respectively. The~p-measure $\mu_1\boxplus\mu_2$ on the~other hand
denotes the~free (additive) convolution of $\mu_1$ and $\mu_2$ introduced by 
Voiculescu~\cite{Vo:1986} for compactly supported p-measures. 
Free convolution was extended by
Maassen~\cite{Ma:1992} to p-measures with finite variance and by Bercovici
and Voiculescu~\cite{BeVo:1993} to the~class $\mathcal M$.
Thus, $\mu_1\boxplus\mu_2$ is the~distribution of $X+Y$,
where $X$ and $Y$ are free random variables with the~distributions $\mu_1$ 
and $\mu_2$, respectively. There are free analogues of multiplicative
convolutions as well; these were first studied in Voiculescu~\cite{Vo:1987}. 

Let $\Bbb C^+\,(\Bbb C^-)$ denote the~open upper (lower) half of
the~complex plane. For $\mu\in\mathcal M$,
define its Cauchy transform by
\begin{equation}\label{2.1}
G_{\mu}(z)=\int\limits_{-\infty}^{\infty}\frac {\mu(dt)}{z-t},
\qquad z\in\Bbb C^+.
\end{equation}

Following Maassen~\cite{Ma:1992} and Bercovici and 
Voiculescu~\cite{BeVo:1993}, we shall consider in the~following
the~ {\it reciprocal Cauchy transform}
\begin{equation}\label{2.2}
F_{\mu}(z)=\frac 1{G_{\mu}(z)}.
\end{equation}
The~corresponding class of reciprocal Cauchy
transforms of all $\mu\in\mathcal M$ will be denoted by $\mathcal F$.
This class admits a~simple description. Recall that the~Nevanlinna
class $\mathcal N$ is the~class  of analytic functions 
$F:\Bbb C^+\to\Bbb C^+\cup\Bbb R$. The~class $\mathcal F$ is
the~subclass of Nevanlinna functions $F_{\mu}$ for which 
$F_{\mu}(z)/z\to 1$ as $z\to \infty$ nontangentially to 
$\infty$ (i.e., such that $\Re z/\Im z$ stays bounded), 
and this implies that $F_{\mu}$ has certain invertibility properties.
(See Akhiezer and Glazman~\cite{AkhG:1963}, 
Akhiezer~\cite {Akh:1965}, Berezanskii~\cite{Ber:1968}). 
More precisely, for two numbers $\alpha>0,\beta>0$ define
$$
\Gamma_{\alpha}=\{z=x+iy\in\Bbb C^+:|x|<\alpha y\}\quad\text{and}\quad
\Gamma_{\alpha,\beta}=\{z=x+iy\in\Gamma_{\alpha}:y>\beta\}.
$$
Then for every $\alpha>0$ there exists $\beta=\beta(\mu,\alpha)$ 
such that $F_{\mu}$ has a~left inverse $F_{\mu}^{(-1)}$ defined on 
$\Gamma_{\alpha,\beta}$.
The~function $\phi_{\mu}(z)=F_{\mu}^{(-1)}(z)-z$ is called
the~Voiculescu transform of $\mu$. It is not hard to show that
$\phi_{\mu}(z)$ is an~analytic function on $\Gamma_{\alpha,\beta}$ and 
$\Im \phi_{\mu}(z)\le 0$ for $z\in \Gamma_{\alpha,\beta}$, where 
$\phi_{\mu}$ is defined. Furthermore, note that $\phi_{\mu}(z)=o(z)$
as $|z|\to\infty$, $z\in\Gamma_{\alpha}$.

Based on alternative definition of free convolution
developed in Chistyakov and G\"otze \cite {ChG:2005},
we define the~free convolution $\mu_1\boxplus\mu_2$ of p-measures
$\mu_1$ and $\mu_2$ as follows. Let $F_{\mu_1}(z)$
and $F_{\mu_2}(z)$ denote their reciprocal Cauchy transforms respectively. 
We shall define the~free convolution $\mu_1\boxplus\mu_2$,
using $F_{\mu_1}(z)$ and $F_{\mu_2}(z)$ only. It was proved
in Chistyakov and G\"otze~\cite {ChG:2005} that
there exist unique functions $Z_1(z)$ and $Z_2(z)$ in the~class 
$\mathcal F$ such that, for $z\in\Bbb C^+$, 
\begin{equation}\label{2.3}
z=Z_1(z)+Z_2(z)-F_{\mu_1}(Z_1(z))\quad\text{and}\quad
F_{\mu_1}(Z_1(z))=F_{\mu_2}(Z_2(z)). 
\end{equation}
The~function $F_{\mu_1}(Z_1(z))$ belongs again to the~class $\mathcal F$ and 
hence by Remark~\ref{3.0} (see Section~3) there exists 
a~p-measure $\mu$ such that
$F_{\mu_1}(Z_1(z)) =F_{\mu}(z)$, where $F_{\mu}(z)=1/G_{\mu}(z)$ and 
$G_{\mu}(z)$ is the~Cauchy transform as in (\ref{2.1}). 
We define $\mu_1\boxplus\mu_2:=\mu$. 
The~measure $\mu$ depends on $\mu_1$ and $\mu_2$ only.

On the~domain $\Gamma_{\alpha,\beta}$, where the~functions $\phi_{\mu_1}(z)$
and $\phi_{\mu_2}(z)$ are defined, we have
\begin{equation}\label{2.4}
\phi_{\mu_1\boxplus\mu_2}(z)=\phi_{\mu_1}(z)+\phi_{\mu_2}(z).
\end{equation} 
This relation  for the~distribution $\mu_1\boxplus\mu_2$ of $X+Y$,
where $X$ and $Y$ are free random variables, is due to
Voiculescu~\cite{Vo:1986}
for the case of  compactly supported p-measures.
The~result was extended by Maassen~\cite{Ma:1992} to p-measures 
with finite variance; the~general case was proved by 
Bercovici and Voiculescu~\cite{BeVo:1993}.
Note that Voiculescu and  Bercovici's definition uses 
the~operator context for the~definition of $\mu_1\boxplus\mu_2$, whereas
Maassen's ~approach is closest to our analytical definition for
the additive free convolution of arbitrary p-measures.
Note that this approach extends as well to the~case 
of multiplicative free convolutions (see \cite {ChG:2005}).
By (\ref{2.4}) it follows  that our definition of $\mu_1\boxplus\mu_2$
coincides with that of~Voiculesku and Bercovici as well as
Maassen's definition.

There is a~notion of infinitely divisible p-measures for additive 
free convolution. As in the~classical case, a~p-measure $\mu$ is 
$\boxplus$-infinitely divisible if, for every natural number $n$,
$\mu$ can be written as $\mu=\nu_n\boxplus\nu_n\boxplus\dots
\boxplus\nu_n$ ($n$ times) with $\nu_n\in\mathcal M$. 
Such $\boxplus$-infinitely divisible p-measures were characterized
by Voiculescu~\cite{Vo:1986} for compactly supported measures.
The $\boxplus$-infinitely divisible p-measures with finite variance were 
studied in Maassen~\cite{Ma:1992} and Bercovichi and 
Voiculescu~\cite{BeVo:1993} extended these results to the~general case. 
There is an~analogue of the~L\'evy-Khintchine formula,
(see Voiculescu, Dykema, Nica~\cite{Vo:1992}, Bercovici and 
Voiculescu~\cite{BeVo:1992}, Bercovici and Voiculescu~\cite{BeVo:1993}
which states that a~p-measure $\mu$, on $\Bbb R$, is infinitely
divisible if and only if the~function $\phi_{\mu}(z)$ has an~analytic
continuation to $\Bbb C^+$, with values in $\Bbb C^-
\cup\Bbb R$, such that
\begin{equation}\label{2.5}
\lim_{y\to +\infty}\frac{\phi_{\mu}(iy)}y=0.
\end{equation}
By the~Nevanlinna representation for such function (see Section~3), 
we know that there exist a~real number $\alpha$, and a~finite nonnegative
measure $\nu$, on $\Bbb R$, such that
\begin{equation}\label{2.6}
\phi_{\mu}(z)=\alpha+\int\limits_{\Bbb R}\frac{1+uz}{z-u}\,\nu(du),
\quad z\in\Bbb C^+.
\end{equation}
Since there is a~one-to-one correspondence between functions
$\phi_{\mu}(z)$ and pairs $(\alpha,\nu)$,
we shall write $\phi_{\mu}=(\alpha,\nu)$.  

Formula~(\ref{2.6}) is an~analogue of the~well-known L\'evy-Khintchine
formula for characteristic functions $\varphi(t;\mu):=\int_{\Bbb R}
e^{itu}\,\mu(du),\,t\in\Bbb R$, of $*$-infinitely
divisible measures $\mu\in\mathcal M$. A~measure $\mu\in\mathcal M$ is
$*$-infinitely divisible if and only if there exist a~finite 
nonnegative Borel measure $\nu$ on $\Bbb R$, and a~real number
$\alpha$ such that
\begin{equation}\label{2.7}
\varphi(t;\mu)=\exp\{f_{\mu}(t)\}:=\exp\Big\{i\alpha t
+\int\limits_{\Bbb R}\Big(e^{itu}-1
-\frac {itu}{1+u^2}\Big)\frac{1+u^2}{u^2}\,\nu(du)\Big\},
\quad t\in\Bbb R,
\end{equation}
where $(e^{itu}-1-itu/(1+u^2))(1+u^2)/u^2$ is defined as
$-t^2/2$ when $u=0$.
Since there is again a~one-to-one correspondence between functions
$f_{\mu}(t)$ and pairs $(\alpha,\nu)$, we shall write $f_{\mu}=(\alpha,\nu)$.  

Bercovici and Pata~\cite{BeP:1999} determined the~distributional
behavior of sums of free {\it identically} distributed 
infinitesimal random variables. More precisely , they showed that, 
given a~sequence $\mu_n$ of p-measures, and an~increasing sequence 
$k_n$ of positive integers, the~free convolution product of $k_n$ measures 
identical to $\mu_n$ converges weakly to a~free infinitely 
divisible distribution
if and only if the corresponding classical convolution product
converges weakly to a~classical infinitely divisible distribution.
Moreover, the~correspondence between the~classical and free limits
can be described explicitly. 

In the~classical case the~precise formulation of the~limit
problem is as follows:

Let $\{\mu_{nk}:\,n\ge 1,\,1\le k\le k_n\}$ be an~array of measures
in $\mathcal M$ such that
\begin{equation}\label{2.8}
\lim_{n\to\infty}\max_{1\le k\le k_n}\mu_{nk}(\{u:|u|>\varepsilon\})=0
\end{equation}
for every $\varepsilon>0$, and $\{a_n:\,n\ge 1\}$ a~sequence of real
numbers. Such triangular schemes of~measures $\mu_{nk}$ are called 
infinitesimal.
The basic limit  problem arising in this context is:
\begin{enumerate}
\item  Find all $\mu\in\mathcal M$ such that $\mu^{(n)}=\delta_{-a_n}*\mu_{n1}
*\mu_{n2}*\dots *\mu_{nk_n}$ converges to $\mu$ in the~weak topology.
\item  Find conditions such that  $\mu^{(n)}$ converges to a given $\mu$.
\end{enumerate}



The~complete solution of this problem has been obtained by the~efforts of
Kolmogorov, P. L\'evy, Feller, de Finetti, Bawly, Khintchine,
Marcinkewicz, Gnedenko, and Doblin.

The~limit problem in free probability theory has the~same form
for the~p-measures $\mu^{(n)}=\delta_{-a_n}\boxplus\mu_{n1}
\boxplus\mu_{n2}\boxplus\dots\boxplus\mu_{nk_n}$.
In the~sequel we denote by $\widehat{\mu}_{nk}$ p-measures such that 
$\widehat{\mu}_{nk}((-\infty,u)):=
\mu_{nk}(-\infty,u+a_{nk}))$, where $a_{nk}:=\int_{(-\tau,\tau)}
u\,\mu_{nk}(du)$ with finite $\tau>0$ which is arbitrary, but fixed.

We provide the~complete solution of this limit problem for 
free random variables.
For the~cla\-ssical case see in Gnedenko and Kolmogorov~\cite{GKo:1968}, 
Ch.~4 and Lo\`eve~\cite{Lo:1963}, \S 22.
\begin{theorem}\label{2.1,1}
Let $\mu_{nk}$ be a triangular scheme of infinitesimal probability measures.
Then we have
\begin{enumerate}
\item  The~family of limit measures of sequences 
$\mu^{(n)}=\delta_{-a_n}\boxplus\mu_{n1}\boxplus\mu_{n2}
\boxplus\dots\boxplus\mu_{nk_n}$
coincides with the~family of $\boxplus$-infinitely divisible 
measures.

\item There exist constants $a_n$ such that the~sequence 
$\mu^{(n)}=\delta_{-a_n}\boxplus\mu_{n1}\boxplus\mu_{n2}
\boxplus\dots\boxplus\mu_{nk_n}$ converges weakly if, and only if, 
$\nu_n$ converges weakly to some finite nonnegative measure $\nu$, where 
$\nu_n$, for any Borel set $S$,
$$
\nu_n(S):=\sum_{k=1}^{k_n}\int\limits_{S}\frac{u^2}{1+u^2}\,
\widehat{\mu}_{nk}(du).
$$
Then all admissible $a_n$ are of the~form $a_n=\alpha_n-\alpha+o(1)$,
where $\alpha$ is an~arbitrary finite number and
$$
\alpha_n=\sum_{k=1}^{k_n}\Big(a_{nk}+\int\limits_{\Bbb R}\frac{u}
{1+u^2}\,\widehat{\mu}_{nk}(du)\Big).
$$
Furthermore, all possible limit measures $\mu\in\mathcal M$ have 
a~Voiculescu transform of type $\phi_{\mu}=(\alpha,\nu)$.
\end{enumerate}
\end{theorem}

Note that the~first statement of the~theorem is due to
Bercovici and Pata~\cite{BeP:2000}. Another proof of this statement,
based on the~theory of Delphic semigroups, has been given by Chistyakov
and G\"otze (\cite{ChG:2005}).
We see that this result is an~obvious 
consequence of the~second statement of the~theorem. 

Comparing the~formulations of the~second statement of Theorem~\ref{2.1,1}
and of the~second statement of the~classical Limit Theorem 
(see Lo\'eve~\cite{Lo:1963}, p. 310), we see that these formulations coincide
for $(\Cal M,\boxplus)$ and $(\Cal M,*)$. 
Therefore the~following result holds, which for the~case of identical 
measures $\mu_{nj},\,j=1,\dots,k_n$, is known as Bercovici-Pata 
bijection~\cite{BeP:1999}.
\begin{theorem}\label{2.2,1}
Let $\mu_{nk}$ be a triangular scheme of infinitesimal probability measures.
There exist constants $a_n$ such that the~sequence 
$\delta_{a_n}\boxplus\mu_{n1}\boxplus\mu_{n2}
\boxplus\dots\boxplus\mu_{nk_n}$ converges weakly to $\mu^{\boxplus}
\in\mathcal M$ such that $\phi_{\mu^{\boxplus}}=(\alpha,\nu)$
if and only if the~sequence $\delta_{a_n}*\mu_{n1}*\mu_{n2}
*\dots*\mu_{nk_n}$ converges weakly to $\mu^*\in\mathcal M$ such that 
$f_{\mu^*}=(\alpha,\nu)$.
\end{theorem}

Let $\mu\in\Cal M$. Denote $\mu^{k*}:=\mu*\dots*\mu$ ($k$ times) and
$\mu^{k\boxplus}:=\mu\boxplus\dots\boxplus\mu$ ($k$ times).
Theorem~\ref{2.2,1} in the {\it identical} case
$\mu_{n1}=\dots=\mu_{nk_n}$ has the~following form.
\begin{corollary}\label{2.2,1a}
Let $\mu_n$ be a~sequence of probability measures.
The~sequence $\mu_n^{k_n\boxplus}$ converges weakly to $\mu^{\boxplus}
\in\mathcal M$ such that $\phi_{\mu^{\boxplus}}=(\alpha,\nu)$
if and only if the~sequence $\mu_n^{k_n*}$ converges weakly 
to $\mu^*\in\mathcal M$ such that $f_{\mu^*}=(\alpha,\nu)$.
\end{corollary}

Bercovici and Pata~\cite{BeP:1999} characterized stable laws and domains
of attraction in free probability theory for the~case of identical
p-measures $\mu_{nj}$ and established the~socalled 
Bercovici-Pata bijection between infinitely divisible limits
in $(\Cal M,*)$ and $(\Cal M,\boxplus)$. In particular they proved
Corollary~\ref{2.2,1a}. Our approach allow us to study
the~case of {\it nonidentical} p-measures $\mu_{nj}$ as well
and to obtain the~results about limiting stable laws.

By Theorem~\ref{2.2,1}, all results concerning the~convergence
of distribution functions of free sums can be reduced to the~corresponding
classical results.
In particular one obtains a~criterion for the~semicircle convergence
(the~case when $\phi_{\mu^{\boxplus}}=(\alpha,\delta_0)$,
a~criterion for the~Marchenko-Pastur convergence  
($\phi_{\mu^{\boxplus}}=(\alpha,\lambda\delta_b),\,\lambda>0, b\ne 0$),
as well as
the~degenerate convergence criterion ($\phi_{\mu^{\boxplus}}=(\alpha,\nu=0)$)
for additive free convolution.
These results generalize the~corresponding results of 
Voiculescu \cite{Vo:1985}, Bercovici and Voiculescu~\cite{BeVo:1995}, 
Maassen~\cite{Ma:1992}, Pata~\cite{P:1996}, Bercovici and Pata~\cite{BeP:1996},
and of Lindsay and Pata~\cite{LiP:1997} to the non-identically
distributed case.

Our analytical approach to the~definition of the~additive free convolution
allows us to give explicit estimates for the~rate of convergence of 
distribution functions of free sums. We shall demonstrate this by proving
a~semicircle approximation theorem (an~analogue
of the~Berry-Esseen inequality (see \cite{Lo:1963}, p. 288), 
and quantitative version of the law of~large numbers, i.e.,
including estimates of convergence.


To formulate the~corresponding results we need the~following notation.
Let $\mu$ be a~p-measure. Define $m_k(\mu):=\int_{\Bbb R}
u^k\,\mu(du)$ and $\beta_k(\mu):=\int_{\Bbb R}|u|^k\,\mu(du)$, where
$k=0,1,\dots$. We denote by $\mu_w$ the~semicircle p-measure, i.e., 
the~measure with the~density $\frac 1{2\pi}
\sqrt{(4-x^2)_+}, \,x\in\Bbb R$, where $a_+:=\max\{a,0\}$ 
for $a\in\Bbb R$. 

Denote by $\Delta(\mu,\nu)$ the~Kolmogorov
distance between the~p-measures $\mu$ and $\nu$, i.e.,
$$
\Delta(\mu,\nu):=\sup_{x\in\Bbb R}|\mu((-\infty,x))-\nu((-\infty,x))|,
$$
and by $L(\mu,\nu)$ the~L\'evy distance between these measures, i.e.,
$$
L(\mu,\nu):=\inf\{h:\mu((-\infty,x-h))-h\le\nu((-\infty,x))\le
\mu((-\infty,x+h))+h,\,x\in\Bbb R\}. 
$$
As it is easy to see, $L(\mu,\nu)\le\Delta(\mu,\nu)$. 

Let $\mu$ be a~p-measure such that $m_1(\mu)=0$ and $m_2(\mu)<\infty$. 
Denote $\mu_n((-\infty,x)):=\mu((-\infty,x\sqrt {m_2(\mu)n})),\,x\in\Bbb R$.

The~following theorem is an~analog of the~well-known Berry-Esseen 
inequality (see \cite{Lo:1963}, p. 288) for the~case of identically
distributed free random variables assuming that the~moment
condition $m_4(\mu)<\infty$ holds.
\begin{theorem}\label{2.3,1}
Let $\mu$ be a~p-measure such that $m_1(\mu)=0$ and 
$m_2(\mu)=1$.  

If $m_4(\mu)<\infty$, there exists an~absolute
constant $c>0$ such that
\begin{equation}\label{2.9,a}
\Delta(\mu_n^{n\boxplus},\mu_w)
\le c\frac {|m_3(\mu)|+(m_4(\mu))^{1/2}}{\sqrt n}.
\end{equation} 
\end{theorem}

The~following proposition shows that estimate (\ref{2.9,a}) is sharp.
\begin{proposition}\label{2.4,1}
Let $\mu$ be a~p-measure such that $\mu(\{-\sqrt{p/q}\})=q$
and $\mu(\{\sqrt{q/p}\})=p$, where $0<p<1,q=1-p$, and $p-q\ne 0$. Then
$$
\Delta(\mu_n^{n\boxplus},\mu_w)\ge L(\mu_n^{n\boxplus},\mu_w)
\ge\frac{c(p)}{\sqrt n},
$$
where $c(p)$ is a~positive constant, depending on $p$ only.
\end{proposition}

Now we shall consider the~case of nonidentically free random variables.
Let $\{\mu_j\}_{j=1}^{\infty}$ be a~sequence of measures in $\mathcal M$
such that $m_1(\mu_j)=0$ and $\beta_3(\mu_j)<\infty$ for all $j=1,\dots$.
Denote
$$
B_n^2=\sum_{k=1}^n m_2(\mu_k),\quad A_n:=\sum_{k=1}^n \beta_3(\mu_k),
\quad L_n:=\frac{A_n}{B_n^3}.
$$
Write $\mu_{nk}((-\infty,x)):=\mu_k((-\infty,B_nx),\,
x\in\Bbb R,\,k=1,\dots,n$,  and $\mu^{(n)}:=\mu_{n1}\boxplus\dots
\boxplus \mu_{nn}$ as well.

\begin{theorem}\label{2.5,1}
There exists an~absolute constant $c>0$ such that
\begin{equation}\label{2.10,b}
\Delta(\mu^{(n)},\mu_w)\le c L_n^{1/2},\quad n=1,\dots.
\end{equation}
\end{theorem}

Finally we shall formulate the~classical degenerate convergence criterion
for additive free convolution with an~estimate of the~convergence.

Let $\{\mu_j\}_{j=1}^{\infty}$ be a~sequence of measures in $\mathcal M$
and let $\mu_{nk}((-\infty,x)):=\mu_k((-\infty,nx))$, $x\in\Bbb R$,
for $k=1,\dots,n$. Denote
$\mu^{(n)}:=\mu_{n1}\boxplus\dots\boxplus\mu_{nn}$.
\begin{theorem}\label{2.6,1}
In order that 
\begin{equation}\label{2.11,a}
L(\mu^{(n)},\delta_0)\to 0
\end{equation} 
as $n\to\infty$ it is necessary and sufficient that, for $n\to\infty$,
\begin{align}
\eta_{n1}&:=\sum_{k=1}^n\int\limits_{\{|x|\ge n\}}\,
\mu_k\to 0,\label {2.11,b}\\
\eta_{n2}&:=\frac 1{n}\sum_{k=1}^n\int\limits_{(-n,n)}
\,x\,\mu_k(dx)\to 0.\label {2.11,cc}\\
\eta_{n3}&:=\frac 1{n^2}\sum_{k=1}^n\Big\{\int\limits_{(-n,n)}
\,x^2\,\mu_k(dx)-\Big(\int\limits_{(-n,n)}\,x\,\mu_k(dx)\Big)^2\Big\}
\to 0,\label {2.11,c}
\end{align}\
In addition, for some absolute positive constant $c$ 
\begin{equation}\label{2.11,d}
L(\mu^{(n)},\delta_0)\le c\Big((\eta_{n1}+\eta_{n3})^{1/6}+\eta_{n2}\big),
\quad n=1,\dots.
\end{equation}
\end{theorem}
 
Note that the~statement of this theorem without the~quantitave bound 
(\ref{2.11,d}) is a simple consequence of Theorem~\ref{2.2,1} and 
the~classical degenerate
criterion (see \cite{Lo:1963}, p. 318). Therefore we need to prove  
(\ref{2.11,d}) only.

Finally we shall describe L\'evy's class $\mathcal L_{\boxplus}$ 
of limit laws of normed sums obeying the~infinitesimal condition
for the~case of free summands.
Let $\mu_1,\mu_2,\dots$ be a~sequence of measures in $\mathcal M$ and
sequences of real numbers $\{a_n\}$ and $\{b_n>0\}$. Denote
by $\mu_{nk}:\,n\ge 1,\,1\le k\le n$, the~measures such that
$\mu_{nk}(S):=\mu_k(b_nS)$ for every Borel set $S$. 

Consider again
the~sequence of measures $\{\mu^{(n)}:=\delta_{-a_n}\boxplus
\mu_{n1}\boxplus\mu_{n2}\boxplus\dots\boxplus\mu_{nn}\}$.

As in the~classical case the~following problems arise.
\begin{enumerate}
\item Given a~sequence $\{\mu_n\}$ of measures in $\mathcal M$, find whether
there exist sequences $\{a_n\}$ and $\{b_n>0\}$ such that the~$\mu_{nk},
\,n\ge 1,\,k=1,\dots,n$, are infinitesimal and $\mu^{(n)}\to \mu$
weakly as $n\to\infty$, where $\mu$ is an~infinitely divisible probability
distribution such that $\phi_{\mu}=(\alpha,\nu)$. If such
sequences exist, then characterize them.

\item Characterize the~family $\mathcal L_{\boxplus}$; in other words, 
characterize those functions $\phi_{\mu}(z)$  and the~corresponding 
measures $\nu$ which represent limit measures of $\mu^{(n)}$. 
\end{enumerate}
It is convenient to exlude degenerate limit distributions 
from our consideration. 
   
In the~first step we give conditions for convergence of $\{\mu^{(n)}\}$.

Let $\widehat{\mu}_{nk},\,n\ge 1,\,k=1,\dots,n$, are p-measures such that
\begin{equation}\label{2.9}
\widehat{\mu}_{nk}((-\infty,u)):=\mu_k((-\infty,b_nu+b_{nk})),
\end{equation}
where $b_{nk}:=\int_{(-b_n,b_n)}x\,\mu_k(dx)$. By $\bar{\mu}_k$ we denote
p-measures such that $\bar{\mu}_k(S)=\mu_k(-S)$
for any Borel set $S$ and by $\mu_k^s$ we denote the~measures
$\mu_k\boxplus\bar{\mu}_k$. Define the~measures
$\nu_n,\,n=1,\dots$, in the~following way. For any Borel set $S$ put
\begin{equation}\label{2.10}
\nu_n(S):=\sum_{k=1}^n\int\limits_{S}\frac {u^2}{1+u^2}\,
\widehat{\mu}_{nk}(du).
\end{equation}

Let us prove the~following norming theorem.
\begin{theorem}\label{2.7,1}
There exist sequence  $b_n>0$ such that 
the~measures $\mu^{(n)}$ converge weakly as $n\to\infty$ to a~nondegenerate
p-measure $\mu$ for suitable $a_n$ and  $\mu_{nk},\,k=1,\dots,n$, 
are infinitesimal, if, and only if, there
exists a~finite nonnegative mesure $\nu$ such that, upon setting in 
$(\ref{2.9})$, $b_n=b_n'>0$ determined by
\begin{equation}\label{2.10a}
\frac 12\sum_{k=1}^n\int\limits_{\Bbb R}\frac{u^2}{(b_n')^2+u^2}\,
\mu_k^s(du)=\nu(\Bbb R),
\end{equation}
we have 
\begin{equation}\label{2.11}
\max_{k=1,\dots,n}\int\limits_{\Bbb R}\frac{u^2}{(b_n')^2+u^2}\,
\mu_k(du)\to 0
\end{equation}
and
\begin{equation}\label{2.12}
\nu_n\to\nu\quad\text{weakly as}\quad n\to\infty.
\end{equation}

The~Voiculescu transform of the~measure $\mu$ has the~form 
$\phi_{\mu}=(\alpha,\nu)$.
\end{theorem}

Comparing this theorem with the~classical result for cumulative
sums (see Gnedenko, Kolmogorov~(\cite{GKo:1968}), \S 31, 
Lo\'eve~(\cite{Lo:1963}), \S 23), similar as in Theorem~\ref{2.2,1}
we obtain the~Bercovi\-ci-Pata bijection for the~case of infinitesimal 
measures $\mu_{nj}$, which are rescaled versions of the measures
$\mu_j$.
\begin{theorem}\label{2.8,1}
There exist constants $a_n$ and $b_n>0$  such that $\mu_{nk},\,
k=1,\dots,n$, are infinitesimal and the~sequence 
$\delta_{a_n}\boxplus\mu_{n1}\boxplus\mu_{n2}
\boxplus\dots\boxplus\mu_{nn}$ converges weakly to $\mu^{\boxplus}
\in\mathcal M$ such that $\phi_{\mu^{\boxplus}}=(\alpha,\nu)$
if and only if the~sequence $\delta_{a_n}*\mu_{n1}*\mu_{n2}
*\dots*\mu_{nn}$ converges weakly to $\mu^*\in\mathcal M$ such that 
$f_{\mu^*}=(\alpha,\nu)$.
\end{theorem}

Using the~classical results about the~class $\mathcal L_*$ (see 
Gnedenko, Kolmogorov~\cite{GKo:1968}, \S 30, Lo\'eve~\cite{Lo:1963}, 
\S 23), we obtain from Theorem~\ref{2.8,1} the~canonical representation
of the~measures of the~class $\mathcal L_{\boxplus}$.

\begin{theorem}\label{2.9,1}
In order that $\mu\in\mathcal M$ belong to the~class $\mathcal L_{\boxplus}$
it is necessary and sufficient that Voiculescu's transform 
of the~measure $\mu$ has the~form $\phi_{\mu}=(\alpha,\nu)$,
where on $(-\infty,0)$ and $(0,\infty)$ the~left and right
derivatives of the~function $\nu(u):=\nu((-\infty,u)),\,
u\in\Bbb R$ , denoted indifferently by $\nu'(u)$, exist
and $\frac{1+u^2}u\nu'(u)$ do not increase.
\end{theorem}

The~class $\mathcal L_{\boxplus}$ admits another discription. 
Let $\mu\in\mathcal M$. Then
for any real constant $\gamma\ne 0$, we
denote by $D_{\gamma}\mu$ the~measure on $\Bbb R$ given by
$D_{\gamma}\mu(S)=\mu(\gamma^{-1}S)$ for any Borel set $S$.

\begin{theorem}\label{2.10,1}
In order that $\mu\in\mathcal M$ belong to the~class $\mathcal L_{\boxplus}$,
it is necessary and sufficient that for every $\gamma,\,0<\gamma<1$,
$\mu=D_{\gamma}\mu\boxplus\mu_{\gamma}$, where $\mu_{\gamma}\in\mathcal M$.
\end{theorem}

\begin{remark}\label{2.11,1}
For any $\gamma\in(0,1)$ the~measure $\mu_{\gamma}$ is 
$\boxplus$-infinitely divisible and $\phi_{\mu_{\gamma}}
=(\alpha_{\gamma},\nu_{\gamma})$. Moreover, for any $\gamma\in(0,1)$, 
$\mu=D_{\gamma}\mu*\mu'_{\gamma}$, where $\mu'_{\gamma}$ is
$*$-infinitely divisible and $f_{\mu'_{\gamma}}=(\alpha_{\gamma},
\nu_{\gamma})$.
\end{remark}



Barndorff-Nielsen and Thorbj{\o}rnsen in \cite{BarTho:2002} and
\cite{BarTho:2004} studied the~connection between  
the~classes of infinitely divisible p-measures in 
classical and free probability. In \cite{BarTho:2002a}  
they studied the~property
of self-decomposa\-bility in free probability and, proving that 
such laws are infinitely divisible, studied L\'evy processes
in free probability and construct stochastic integrals with respect
to such processes. Our results allow to extend the~results of 
Biane~\cite{Bi:1998} and of Barndorff-Nielsen and 
Thorbj{\o}rnsen~\cite{BarTho:2002}--\cite{BarTho:2004}. 

\section{Auxilliary results}

We need results about some classes of analytic functions
(see {\cite{Akh:1965}, Section~3, and {\cite{AkhG:1963}},
Section~6, \S 59). 

The~class $\mathcal N$ (Nevanlinna, R.) is the~class of analytic 
functions $f(z):\Bbb C^+\to\{z: \,\Im z\ge 0\}$.
For such functions there is the~integral representation
\begin{equation}\label{3.1}
f(z)=a+bz+\int\limits_{\Bbb R}\frac{1+uz}{u-z}\,\tau(du)=
a+bz+\int\limits_{\Bbb R}\Big(\frac 1{u-z}-\frac u{1+u^2}\Big)(1+u^2)
\,\tau(du),\quad z\in\Bbb C^+,
\end{equation}
where $b\ge 0$, $a\in\Bbb R$, and $\tau$ is a~nonnegative finite
measure. Moreover, $a=\Re f(i)$ and $\tau(\Bbb R)=\Im f(i)-b$.   

A~function $f\in\mathcal N$ admits the~representation
\begin{equation}\label{3.2}
f(z)=\int\limits_{\Bbb R}\frac{\sigma(du)}{u-z},\quad z\in\Bbb C^+,
\end{equation}
where $\sigma$ is a~finite nonnegative measure, if and only if
$\sup_{y\ge 1}|yf(iy)|<\infty$.

\begin{remark}\label{3.0}
Since the~class $\mathcal F$ is the~subclass of Nevalinna functions
$f(z)$ for which $f(z)/z\to 1$ as $z\to\infty$ nontangentially,
we note that every $f\in\mathcal F$ admits representation (\ref{3.1}),
where $b=1$. Moreover $-1/f(z)$ admits representation (\ref{3.2}),
where $\sigma\in\mathcal M$. Note as well that a~function $f\in\Cal F$
satifies the~obvious inequality
\begin{equation}\label{3.2,a}
\Im f(z)\ge \Im z,\quad z\in\Bbb C^+.
\end{equation}
\end{remark}

The~Stieltjes-Perron inversion formula for the~functions $f$ of class 
$\mathcal N$ has the~following form.
Let $\psi(u):=\int_0^u(1+t^2)\,\tau(dt)$. Then
\begin{equation}\label{3.3}
\psi(u_2)-\psi(u_1)=\lim_{\eta\to 0}\frac 1{\pi}\int\limits_{u_1}
^{u_2}\Im f(\xi+i\eta)\,d\xi,
\end{equation}
where $u_1<u_2$ denote two continuity points of the~function $\psi(u)$.

Furthermore, we shall need the~following inequality for the~distance 
between distributions in terms of their Stieltjes transform. 
\begin{lemma}\label{3.1a}
Let $\mu_w$ be the~semicircle measure and $\mu$ be a~p-measure
such that 
\begin{equation}\label{3.3a}
\int\limits_{\Bbb R}|\mu_w((-\infty,x))-\mu((-\infty,x))|\,dx<\infty.
\end{equation}
Then there exists an~absolute constant $c$ such that, for any $0<v<1$,
\begin{align}
\Delta(\mu_w,\mu)
\le c\int\limits_{\Bbb R}|G_{\mu_w}(u+i)&-G_{\mu}(u+i)|\,du+cv\notag\\
&+c\sup_{x\in[-2,2]}|\int\limits_v^1(G_{\mu_w}(x+iu)
-G_{\mu}(x+iu))\,du|,
\notag
\end{align}
where $G_{\mu_w}$ and $G_{\mu}$ are defined in $(\ref{2.1})$.
\end{lemma}
This lemma is a~simple consequence of Corollary~2.3 in G\"otze and 
Tikhomirov \cite{GT:2003}.

Let $\mu_j\in\mathcal M,\,j=1,2$. Recall that $\mu_1\boxplus\mu_2$ 
is defined in Chistyakov and G\"otze \cite{ChG:2005} as follows.

\begin{theorem}\label{3.1b}
Let $F_{\mu_1}(z)$ and $F_{\mu_1}(z)$ be the~reciprocal Cauchy transforms
of the~p-measures $\mu_1$ and $\mu_2$, respectively. Then
there exist unique functions $Z_1(z)$ and $Z_2(z)$ of class 
$\mathcal F$ such that 
$$
z=Z_1(z)+Z_2(z)-F_{\mu_1}(Z_1(z))\quad\text{and}\quad
F_{\mu_1}(Z_1(z))=F_{\mu_2}(Z_2(z)),\quad z\in\Bbb C^+.
$$
\end{theorem}

Consider the~function $F_{\mu_1}(Z_1(z)),\,z\in\Bbb C^+$. It is easy 
to see that this function belongs to the~class $\mathcal F$. Therefore
there exists a~measure $\mu\in\mathcal M$ such that $F_{\mu}(z)=1/G_{\mu}(z)
=F_{\mu_1}(Z_1(z)),z\in \Bbb C^+$. We define $\mu:=\mu_1\boxplus\mu_2$.

As shown in \cite{ChG:2005}, Theorem~\ref{3.1b} admits the~following
consequences.
\begin{corollary}\label{3.2a}
Let $\mu_1,\dots,\mu_n\in\Cal M$.
There exist unique functions $Z_1(z),\dots,Z_n(z)$ of class 
$\mathcal F$ such that, for $z\in\Bbb C^+$,
\begin{equation}\label{3.7}
z=Z_1(z)+\dots+Z_n(z)-(n-1)F_{\mu_1}(Z_1(z)),\quad\text{and}
\quad F_{\mu_1}(Z_1(z))=\dots=F_{\mu_n}(Z_n(z)).
\end{equation}
Moreover, $F_{\mu_1\boxplus\dots\boxplus\mu_n}(z)=F_{\mu_1}(Z_1(z))$
for all $z\in\Bbb C^+$.
\end{corollary}

Let $\mu_1=\mu_2=\dots=\mu_n=\mu$ and write ${\mu_1\boxplus\dots\boxplus\mu_n}
=\mu^{n\boxplus}$.
\begin{corollary}\label{3.3b}
Let $\mu\in\mathcal M$. There exists a~unique function $Z\in\mathcal F$ 
such that
\begin{equation}\label{3.8}
z=nZ(z)-(n-1)F_{\mu}(Z(z)),\quad z\in\Bbb C^+,
\end{equation}
and $F_{\mu^{n\boxplus}}(z)=F_{\mu}(Z(z)),\,z\in\Bbb C^+$.
\end{corollary}

We need the~following auxiliary results of Bercovici and 
Voiculescu~\cite{BeVo:1993}.




\begin{proposition}\label{3.4a}
Let $\big\{\mu_n\big\}_{n=1}^{\infty}$ be a~sequence of p-measures 
on $\Bbb R$. The~fo\-llowing assertions are equivalent.
\begin{enumerate}
\item  The~sequence $\big\{\mu_n\big\}_{n=1}^{\infty}$ converges  
weakly to a~p-measure $\mu$.

\item  There exist $\alpha,\beta>0$ such that the~sequence
$\big\{\phi_{\mu_n}\big\}_{n=1}^{\infty}$ converges uniformly
on compact subsets of $\Gamma_{\alpha,\beta}$ to a~function 
$\phi$, and $\phi_{\mu_n}(iy)=o(y)$ uniformly in $n$ as
$y\to+\infty$. 
\end{enumerate}
Moreover, if $(a)$ and
$(b)$ are satisfied, we have $\phi=\phi_{\mu}$ in 
$\Gamma_{\alpha,\beta}$. 
\end{proposition}

\begin{proposition}\label{3.5a}
Let $\big\{\mu_n\big\}_{n=1}^{\infty}$ and 
$\big\{\nu_n\big\}_{n=1}^{\infty}$ be sequences of p-measures 
on $\Bbb R$ which converge weakly to p-measures 
$\mu$ and $\nu$, respectively. Then $\big\{\mu_n\boxplus\nu_n
\big\}_{n=1}^{\infty}$ converges weakly to the~p-measure
$\mu\boxplus\nu$.
\end{proposition}

We also need the~following two results which are due to Bercovici
and Pata~\cite{BeP:2000}, \cite{BeP:1999}.
\begin{proposition}\label{3.6a}
Let $\alpha,\beta,\varepsilon$ be positive numbers, and let 
$\phi:\Gamma_{\alpha,\beta}\to\Bbb C$
be an~analytic function such that
$$
|\phi(z)|\le\varepsilon |z|,\qquad z\in\Gamma_{\alpha,\beta}.
$$
For every $\alpha'<\alpha$ and $\beta'>\beta$ there exists $k>0$ such that
$$
|\phi'(z)|\le k\varepsilon.
$$
\end{proposition}

\begin{proposition}\label{3.6aa}
For every $\alpha,\beta>0$ there exists $\varepsilon>0$ with 
the~following property. If $\mu\in\mathcal M$ such that $\int
_{\Bbb R}u^2/(1+u^2)\,\mu(du)<\varepsilon$, then $\phi_{\mu}$
is defined on the~region $\Gamma_{\alpha,\beta}$ and $\phi_{\mu}(
\Gamma_{\alpha,\beta})\subset \Bbb C^-\cup \Bbb R$.
\end{proposition}

Let $\mu$ be a~p-measure. Denote by $\bar\mu$ the~measure
defined by $\bar\mu(B)=\mu(-B)$ for any Borel set $B$. Write
$\mu^s:=\mu\boxplus\bar\mu$.

\begin{proposition}\label{3.7a}
A~p-measure $\mu$ is symmetric if and only if the~functions
$G_{\mu}(iy)$ and $F_{\mu}(iy)$ take imaginary values for $y>0$
and the~function $\phi_{\mu}(iy)$ takes imaginary values on the~set
$y\ge y_0>0$, where it is defined.
\end{proposition}
We omit the~proof of this simple proposition.

We obtain, as an~obvious consequences of Proposition~\ref{3.7a}, that
$\mu^s$ is a~symmetric p-measure. In addition,
if $\mu_1$ and $\mu_2$ are symmetric p-measures, then
$\mu_1\boxplus\mu_2$ is a~symmetric p-measure as well.

The following auxiliary result is due to Khintchine (see \cite{GKo:1968},
p. 42).
\begin{proposition}\label{3.8a}
For a~sequence of distribution functions $F_n(x)$ the~relations
$$
F_n(b_nx+a_n)\to F(x)\quad\text{weakly},\quad
F_n(\beta_n x+\alpha_n)\to F(x)\quad\text{weakly}, 
$$
as $n\to\infty$, where $b_n>0,\beta_n>0$, $a_n,\alpha_n$ are real 
constants and $F(x)$ is a~proper distribution function, are
satisfied simultaneously if and only if $(\beta_n/b_n)\to 1$ and
$(a_n-\alpha_n)/b_n\to 0$ as $n\to\infty$.
\end{proposition}

\section{Additive free limit theorem}
In this section we shall prove Theorem~\ref{2.1,1}. In the~sequel
we denote by $c$ positive absolute constants. For some measure $\nu$
and for some parameter $\tau$ we denote by $c(\nu),\,c(\tau)$, and 
$c(\nu,\tau)$ positive constants which only depend on the~measure 
$\nu$, on the~parameter $\tau$, and  on $\nu$ and $\tau$, respectively. 
Before to prove Theorem~\ref{2.1,1}
we establish some properties of the~measures $\{\mu_{nk}:n\ge 1,
1\le k\le k_n\}$, satisfying condition (2.8), and the~corresponding
reciprocal Cauchy transforms $\{F_{\mu_{nk}}(z):n\ge 1$,$1\le k\le k_n\}$.

It is well-known that condition (2.8) is equivalent to the~following 
relation (see Lo\`eve~\cite{Lo:1963}, p.~302)
$$
\max_{k=1,\dots,k_n}\int\limits_{\Bbb R}\frac{u^2}{1+u^2}\,
\mu_{nk}(du)\to 0,\quad n\to\infty.
$$
Recall that $\widehat{\mu}_{nk}((-\infty,u)):=\mu_{nk}((-\infty,u+a_{nk}))$,
where $a_{nk}:=\int_{(-\tau,\tau)}x\,\mu_{nk}(dx)$, $k=1,\dots,k_n$,
with arbitrary $\tau>0$ which is finite and fixed.
Since obviously $\max_{k=1,\dots,k_n}|a_{nk}|\to 0$ as 
$n\to\infty$, we obtain 
\begin{equation}\label{4.1}
\varepsilon_n:=\max_{k=1,\dots,k_n}\varepsilon_{nk}\to 0,
\quad n\to\infty,\quad\text{where}\quad
\varepsilon_{nk}:=\int\limits_{\Bbb R}\frac{u^2}{1+u^2}\,
\widehat{\mu}_{nk}(du).
\end{equation}

By Remark~\ref{3.0}, for every $k=1,\dots,k_n$ the~reciprocal of 
the~Cau\-chy transform of $G_{\widehat{\mu}_{nk}}(z)$, see (\ref{2.1}),
has the form
\begin{equation}\label{4.2}
F_{\widehat{\mu}_{nk}}(z)=b_{nk}+z+\int\limits_{\Bbb R}
\frac{1+uz}{u-z}\,\sigma_{nk}(du),
\end{equation}
where $b_{nk}:=\Re (G_{\widehat{\mu}_{nk}}(i))^{-1}$ and 
$\sigma_{nk}$ is a~nonnegative finite
measure such that $\sigma_{nk}(\Bbb R)
=\Im (G_{\widehat{\mu}_{nk}}(i))^{-1}-1$.  
From (\ref{4.2}) we deduce the~following relation
\begin{equation}\label{4.2a}
-\frac{\Im G_{\widehat{\mu}_{n k}}(iy)}
{|G_{\widehat{\mu}_{n k}}(iy)|^2}=y\Big(1+\int\limits_{\Bbb R}
\frac{1+u^2}{u^2+y^2}\sigma_{nk}(du)\Big),\quad y>0,\,\,k=1,\dots,k_n, 
\end{equation}
which yields
\begin{equation}\label{4.3}
1+\int\limits_{\Bbb R}\frac{1+u^2}{u^2+y^2}\sigma_{nk}(du)
\le -\frac 1{y\Im G_{\widehat{\mu}_{n k}}(iy)},\quad y>0,\,\,k=1,\dots,k_n.
\end{equation}

On the~other hand we see that, for $y\ge 1$,
\begin{equation}\label{4.3a}
-y\Im G_{\widehat{\mu}_{n k}}(iy)=\int\limits_{\Bbb R}\frac{y^2}{u^2+y^2}
\widehat{\mu}_{n k}(du)=1-\int\limits_{\Bbb R}\frac{u^2}{u^2+y^2}
\widehat{\mu}_{n k}(du):=1-\varepsilon_{nk}(y).
\end{equation}
Hence, for sufficiently large $n\ge n_0$ and $k=1,\dots,k_n$, 
we obtain, by (\ref{4.3}) and
(\ref{4.3a}), the~upper bound
\begin{equation}\label{4.4}
\int\limits_{\Bbb R}\frac{1+u^2}{u^2+y^2}\sigma_{nk}(du)\le 
\frac{\varepsilon_{nk}(y)}{1-\varepsilon_{nk}(y)}\le 2\varepsilon_{nk}(y),
\quad y\ge 1.
\end{equation}
It follows from (\ref{4.4}) that, for $n\ge n_0$, 
\begin{equation}\label{4.5}
\sigma_{nk}(\Bbb R)\le 2\varepsilon_{nk}(1)
=2\varepsilon_{nk},\quad k=1,\dots,k_n,
\end{equation}
and $\max_{k=1,\dots,k_n}\sigma_{nk}(\Bbb R)\to 0$ as $n\to\infty$.

Now we deduce the~relation 
\begin{align}
\Re G_{\widehat{\mu}_{nk}}(i)&=\int\limits_{\Bbb R}\frac {u-a_{nk}}
{(u-a_{nk})^2+1}\,\mu_{nk}(du)\notag\\
&=\int\limits_{(-\tau,\tau)}
\frac {u-a_{nk}}{(u-a_{nk})^2+1}\,\mu_{nk}(du)
+\int\limits_{|u|\ge \tau}\frac {u-a_{nk}}{(u-a_{nk})^2+1}\,
\mu_{nk}(du)\notag\\
&=\int\limits_{(-\tau,\tau)}(u-a_{nk})\,\mu_{nk}(du)
-\int\limits_{(-\tau,\tau)}\frac {(u-a_{nk})^3}{(u-a_{nk})^2+1}\,
\mu_{nk}(du)\notag\\
&+\int\limits_{|u|\ge \tau}\frac {u-a_{nk}}
{(u-a_{nk})^2+1}\,\mu_{nk}(du)\notag\\
&=-\int\limits_{(-\tau,\tau)}\frac {(u-a_{nk})^3}{(u-a_{nk})^2+1}\,
\mu_{nk}(du)
+\int\limits_{|u|\ge \tau}\frac {u+a_{nk}(u-a_{nk})^2}{(u-a_{nk})^2+1}\,
\mu_{nk}(du).\notag
\end{align}
Using straighforward estimates we easily have, for sufficiently large
$n\ge n_0$,
\begin{equation}\label{4.6}
|\Re G_{\widehat{\mu}_{nk}}(i)|\le c(\tau)\varepsilon_{nk},\quad k=1,\dots,k_n.
\end{equation}
In view of (\ref{4.1}), (\ref{4.3a}), and (\ref{4.6}), we get, 
for $n\ge n_0$ and $k=1,\dots,k_n$, 
\begin{equation}\label{4.7}
|b_{nk}|\le|\Re G_{\widehat{\mu}_{nk}}(i)|/
(\Im G_{\widehat{\mu}_{nk}}(i))^2
\le c(\tau)\varepsilon_{nk}.
\end{equation}
From (\ref{4.2}), (\ref{4.5}), and (\ref{4.7}) we obtain, for $z\in\Bbb C^+$
and $n\ge n_0,\,k=1,\dots,k_n$,
\begin{align}\label{4.8}
|F_{\widehat{\mu}_{nk}}(z)-z|&\le |b_{nk}|+\int\limits_{\Bbb R}
\frac {\sigma_{nk}(du)}{|u-z|}+
\int\limits_{\Bbb R}\frac{|z||u|}{|u-z|}\,\sigma_{nk}(du)\notag\\
&\le c(\tau)\varepsilon_{nk} \Big(1+\frac{1+|z|^2}{\Im z}\Big)
\le c(\tau)\varepsilon_{nk}Q(z),
\end{align}
where $Q(z):=\frac{1+|z|^2}{\Im z}$.
In addition to (\ref{4.4}) we deduce the~estimate
\begin{align}\label{4.8a}
\Im (F_{\widehat{\mu}_{nk}}&(z)-z)=\Im z\int\limits_{\Bbb R}\frac{1+u^2}
{(u-\Re z)^2+(\Im z)^2}\,\sigma_{nk}(du)\notag\\
&\le 2\Big(\frac{|z|}{\Im z}\Big)^2\Im z
\int\limits_{\Bbb R}\frac{1+u^2}
{(\Im z)^2+u^2}\,\sigma_{nk}(du)
\le 4\Big(\frac{|z|}{\Im z}\Big)^2
\Im z\int\limits_{\Bbb R}\frac{1+u^2}
{(\Im z)^2+u^2}\frac{u^2}{1+u^2}\,\widehat{\mu}_{nk}(du)\notag\\
&\le 4\Big(\frac{|z|}{\Im z}\Big)^2
\Big(\frac 2{\Im z}\int\limits_{[-\sqrt{\Im z},\sqrt{\Im z}]}
\frac{u^2}{1+u^2}\,\widehat{\mu}_{nk}(du)+
\int\limits_{|u|>\sqrt{\Im z}}\frac{u^2}{1+u^2}\,\widehat{\mu}_{nk}(du)
\Big)\Im z\notag\\
&:=4\Big(\frac{|z|}{\Im z}\Big)^2\eta_{nk}(\Im z)\Im z
\end{align}
for $k=1,\dots,k_n$ and $\Im z\ge 1$. Note that, for such $k$ and $\Im z$,
$\eta_{nk}(\Im z)\le 4\varepsilon_{nk}$. 

We conclude from (\ref{4.8}), (\ref{4.8a}) and Rouch\'e's theorem that 
for every $y\ge 1$ there exists
a~neighborhood $|z-iy|\le y/2$ such that the~inverse function
$F_{\widehat{\mu}_{nk}}^{(-1)}(z)$ exists and is analytic in this domain. 
In addition, the~following inequalities hold
\begin{equation}\label{4.9}
|\Re\phi_{\widehat{\mu}_{nk}}(z)|=|\Re(F_{\widehat{\mu}_{nk}}^{(-1)}(z)-z)|
\le c(\tau)\varepsilon_{nk}y,\quad
|\Im\phi_{\widehat{\mu}_{nk}}(z)|=|\Im(F_{\widehat{\mu}_{nk}}^{(-1)}(z)-z)| 
\le c\tilde{\eta}_{nk}(y)y,
\end{equation}
for $|z-iy|\le y/2, \,\,n\ge n_0,\,\,k=1,\dots,k_n$, where 
$\tilde{\eta}_{nk}(y):=\max_{t\in[y/4,2y]}\eta_{nk}(t)$.

\vspace {0,5 cm}

{\bf Proof of Theorem~$\ref{2.1,1}$.} {\it Sufficiency}.
Consider the~measure 
$\widehat{\mu}_n:=\widehat{\mu}_{n1}\boxplus\dots\boxplus\widehat{\mu}_{nk_n}$.
It follows from Corollary~\ref{3.2a} that there exist unique functions
$Z_{n1},\dots,Z_{nk_n}$ of class $\mathcal F$ such that, 
for $ z\in\Bbb C^+$,
\begin{equation}\label{4.10}
Z_{n1}(z)-z=F_{\widehat{\mu}_{n2}}(Z_{n2}(z))-Z_{n2}(z)+\dots
+F_{\widehat{\mu}_{nk_n}}(Z_{nk_n}(z))-Z_{nk_n}(z),
\end{equation}
and
\begin{equation}\label{4.11}
F_{\widehat{\mu}_{n1}}(Z_{n1}(z))=F_{\widehat{\mu}_{n2}}(Z_{n2}(z))=\dots
=F_{\widehat{\mu}_{nk_n}}(Z_{nk_n}(z)).
\end{equation}
Moreover, we have $F_{\widehat{\mu}_n}(z)=F_{\widehat{\mu}_{nk}}(Z_{nk}(z)),
\,z\in\Bbb C^+,\,k=1,\dots,k_n$.
Then, by (\ref{4.9})--(\ref{4.11}), for $|z-iy|\le y/2,\,y\ge 1$,
it follows that
\begin{equation}\label{4.12}
|\Re\phi_{\widehat{\mu}_{n1}\boxplus\dots\boxplus\widehat{\mu}_{nk_n}}(z)|
\le |\Re\phi_{\widehat{\mu}_{n1}}(z)|+\dots+|\Re\phi_{\widehat{\mu}_{nk_n}}(z)|
\le c(\tau)\eta_n y:=c(\tau)\Big(\sum_{k=1}^{k_n}\varepsilon_{nk}\Big)y
\end{equation}
and
\begin{equation}\label{4.12a}
|\Im\phi_{\widehat{\mu}_{n1}\boxplus\dots\boxplus\widehat{\mu}_{nk_n}}(z)|
\le |\Im\phi_{\widehat{\mu}_{n1}}(z)|+\dots+|\Im\phi_{\widehat{\mu}_{nk_n}}(z)|
\le c\eta_n(y)y:=c\Big(\sum_{k=1}^{k_n}\tilde{\eta}_{nk}(y)\Big)y.
\end{equation}
By the~assumptions of the~theorem, we have $\eta_n\le 2\nu(\Bbb R)$
for sufficiently large $n\ge n_0$. In addition, by
(\ref{4.8a}) and the~assumptions of the~theorem, we see that
\begin{equation}\label{4.12b}
\eta_n(y)\le \frac {16\nu(\Bbb R)}y+4\nu(\Bbb R\setminus[-\sqrt y /2,
\sqrt y /2]) 
\end{equation}
for sufficiently large $n\ge n_1(y)$, where $-\sqrt y/2$ and $\sqrt y/2$ 
are continuity points of the~function $\nu((-\infty,x)),\,x\in\Bbb R$. 
In the~sequel we choose $y$ so that $-\sqrt y/2$ and $\sqrt y/2$ are 
continuity points of $\nu((-\infty,x))$.
Since 
$$
\phi_{\widehat{\mu}_{n1}\boxplus\dots\boxplus\widehat{\mu}_{nk_n}}(z)
=(F_{\widehat{\mu}_{nk}}(Z_{nk}))^{(-1)}(z)-z=
Z_{nk}^{(-1)}(F_{\widehat{\mu}_{nk}}^{(-1)}(z))-z,
\,\,k=1,\dots,k_n,
$$ 
for $|z-iy|\le y/2$, we have, by (\ref{4.8}), the~relation
$$
\phi_{\widehat{\mu}_{n1}\boxplus\dots\boxplus\widehat{\mu}_{nk_n}}
(F_{\widehat{\mu}_{nk}}(z))=Z_{nk}^{(-1)}(z)-F_{\widehat{\mu}_{nk}}(z),
\,\,k=1,\dots,k_n,
$$
for $|z-iy|\le y/4$.
Therefore we conclude by (\ref{4.8})--(\ref{4.9}) and 
(\ref{4.12})--(\ref{4.12b}) that the~functions 
$Z_{nk}^{(-1)}(z)$ are analytic in the~disk $|z-iy|<y/4$ and 
\begin{equation}\label{4.13}
|\Re(Z_{nk}^{(-1)}(z)-z)|\le c(\tau)\nu(\Bbb R)Q(y),\quad
|\Im(Z_{nk}^{(-1)}(z)-z)|\le c\Big(\nu(\Bbb R)+y\nu(\Bbb R\setminus
[-\sqrt y/2,\sqrt y/2])\Big),
\end{equation}
for $|z-iy|\le y/4,\,n\ge n_1(y),\,k=1,\dots,k_n$. 
We conclude from (\ref{4.13}) that there exists $y_0=y_0(\nu)\ge 1$ 
such that $Z_{nk}^{(-1)}(z)\in R_{y_0}:=\{z:|\Re z|\le c(\tau)\nu(\Bbb R)y_0,
\,y_0/2\le\Im z\le 3y_0/2\}$ for $|z-iy_0|\le y_0/4,\,n\ge n_1(y_0),\,
k=1,\dots,k_n$. 
Hence there exist points $z_{nk}\in R_{y_0}$ such that $|Z_{nk}(z_{nk})-iy_0|
\le y_0/4$ for $n\ge n_1(y_0),\,k=1,\dots,k_n$. 

The~functions $Z_{nk}$ are of class $\mathcal F$. Therefore
\begin{equation}\label{4.14}
Z_{nk}(z)=d_{nk}+z+\int\limits_{\Bbb R}\frac{1+uz}{u-z}\,
\nu_{nk}(du)=d_{nk}+z+\int\limits_{\Bbb R}\Big(\frac 1{u-z}
-\frac u{1+u^2}\Big)(1+u^2)\,\nu_{nk}(du) 
\end{equation}
for $z\in\Bbb C^+$, where $d_{nk}\in\Bbb R$ and $\nu_{nk}$ 
are finite nonnegative measures.
Since $\Im Z_{nk}(z_{nk})-y_0\le y_0/2$, we have
\begin{equation}\label{4.15a}
c(\nu,\tau)\nu_{nk}(\Bbb R)\le \Im z_{nk}\int\limits_{\Bbb R}\frac{1+u^2}
{(u-\Re z_{nk})^2+(\Im z_{nk})^2}\,\nu_{nk}(du)\le \Im Z_{nk}(z_{nk}) 
\le \frac{3y_0}2.
\end{equation}
It is easy to see from (\ref{4.14}) and (\ref{4.15a}) 
that $|Z_{nk}(z_{nk})-d_{nk}|\le c(\nu,\tau)$. Hence, using the~bound
$|Z_{nk}(z_{nk})|\le 3y_0/2$, we conclude that $|d_{nk}|\le c(\nu,\tau)+3y_0/2$.
Hence we have 
\begin{equation}\label{4.15}
|d_{nk}|\le c(\nu,\tau)\quad\text{and}\quad 
\nu_{nk}(\Bbb R)\le c(\nu,\tau),\quad n\ge n_1(y_0),\,\,k=1,\dots,k_n.
\end{equation}
In the~sequel we assume that $ n\ge n_1(y_0)+n_0$.
As in (\ref{4.8}) we obtain,
for $z\in\Bbb C^+$ and $k=1,\dots,k_n$,
\begin{equation}\label{4.16}
|Z_{nk}(z)-z|\le c(\nu,\tau)Q(z).
\end{equation}
Using (\ref{4.16}) and the~inequality $\Im Z_{nk}(z)\ge \Im z,z\in\Bbb C^+$,
see (\ref{3.2,a}), we deduce 
\begin{equation}\label{4.16a}
Q(Z_{nk}(z))=\frac{1+|Z_{nk}(z)|^2}{\Im Z_{nk}(z)}
\le c(\nu,\tau)\frac 1{\Im z}Q^2(z)
\end{equation}
for $z\in\Bbb C^+$ and $k=1,\dots,k_n$. Therefore we obtain from (\ref{4.8})
\begin{equation}\label{4.17}
|F_{\mu_{nk}}(Z_{nk}(z))-Z_{nk}(z)|\le c(\tau)\varepsilon_{nk}Q(Z_{nk}(z)) 
\le c(\nu,\tau)\varepsilon_{nk}\frac 1{\Im z}Q^2(z)
\end{equation}
for $z\in\Bbb C^+$ and $k=1,\dots,k_n$. 
Let us return to the~relation (\ref{4.11}). In view of (\ref{4.17}), we have,
for $z\in\Bbb C^+$ and $k=1,\dots,k_n$,
\begin{align}\label{4.18}
|Z_{n1}(z)-Z_{nk}(z)|&\le |F_{\widehat{\mu}_{n1}}(Z_{n1}(z))-Z_{n1}(z)|+ 
|F_{\widehat{\mu}_{nk}}(Z_{nk}(z))-Z_{nk}(z)|\\
&\le  c(\nu,\tau)\varepsilon_n\frac 1{\Im z}Q^2(z).
\notag
\end{align}
On the~other hand, in view of (\ref{4.2}) and (\ref{4.5}), we conclude
\begin{align}
&|(F_{\widehat{\mu}_{nk}}(Z_{nk}(z))-Z_{nk}(z))
-(F_{\widehat{\mu}_{nk}}(Z_{n1}(z))-Z_{n1}(z))|\notag\\
&\le\int\limits_{\Bbb R} \frac{|Z_{nk}(z)-Z_{n1}(z)|(1+u^2)\,
\sigma_{nk}(du)}{\sqrt{(u-\Re Z_{nk}(z))^2+(\Im Z_{nk}(z))^2}
\sqrt{(u-\Re Z_{n1}(z))^2+(\Im Z_{n1}(z))^2}}
\notag\\
&\le c\varepsilon_{nk}|Z_{nk}(z)-Z_{n1}(z)|
\frac {(1+|Z_{n1}(z)|)(1+|Z_{nk}(z)|)}{\Im Z_{n1}(z)\Im Z_{nk}(z)} \notag
\end{align}
for $z\in\Bbb C^+$ and $k=1,\dots,k_n$. Thus, taking into 
account (\ref{4.16}), (\ref{4.18}), and the~inequality $\Im Z_{nk}(z)\ge \Im z,\,z\in\Bbb C^+$,
we have, for the~same $z$ and $k$
as above,
\begin{equation}\label{4.19}
|(F_{\widehat{\mu}_{nk}}(Z_{nk}(z))-Z_{nk}(z))
-(F_{\widehat{\mu}_{nk}}(Z_{n1}(z))-Z_{n1}(z))|
\le c(\nu,\tau)\varepsilon_{nk}\varepsilon_n\frac 1{(\Im z)^3}Q^4(z).
\end{equation}

Consider the~functions
$$
f_{nk}(z):=z^2\Big(G_{\widehat{\mu}_{nk}}(z)-\frac 1z\Big)=\gamma_{nk}
+\int\limits_{\Bbb R}\frac{1+uz}{z-u}\,\rho_{nk}(du),
\quad z\in\Bbb C^+,\,\,k=1,\dots,k_n,
$$
where
$$
\gamma_{nk}:=\int\limits_{\Bbb R}\frac u{1+u^2}\,\widehat{\mu}_{nk}(du)
\quad\text{and}\quad
\rho_{nk}(du):=\frac{u^2}{1+u^2}\,\widehat{\mu}_{nk}(du).
$$
By (\ref{4.6}), the~constants $\gamma_{nk}$ admit the~estimates 
$|\gamma_{nk}|\le c(\tau)\varepsilon_{nk}$ for $k=1,\dots,k_n$. 
Hence $\gamma_n:=\sum_{k=1}^{k_n}\gamma_{nk}$ satisfies the~inequality
\begin{equation}\label{4.19a}
|\gamma_n|\le c(\nu,\tau), \quad n=1,\dots.
\end{equation} 
As in (\ref{4.8}), we conclude that 
\begin{equation}\label{4.20}
|f_{nk}(z)|\le c(\tau)\varepsilon_{nk}Q(z),\qquad
z\in\Bbb C^+,\,\,k=1,\dots,k_n.
\end{equation}
We have, for $z\in\Bbb C^+$ and $k=1,\dots,k_n$,
\begin{equation}\label{4.21}
F_{\widehat{\mu}_{nk}}(z)=\frac{z^2}{z+f_{nk}(z)}=
z-f_{nk}(z)+\theta_{nk}(z),
\end{equation}
where 
$$
\theta_{nk}(z)=\frac{f_{nk}^2(z)}{z+f_{nk}(z)}=f_{nk}^2(z)
F_{\widehat{\mu}_{nk}}(z)z^{-2}.
$$
Hence, by (\ref{4.8}) and (\ref{4.20}), 
we conclude, for those $z$, $k$,
\begin{equation}\label{4.22}
|\theta_{nk}(z)|\le c(\tau)\varepsilon_{nk}^2
Q^2(z)\Big(|z|+\varepsilon_{nk}Q(z)\Big)\frac 1{|z|^2}.
\end{equation}
We see from (\ref{4.16}), (\ref{4.16a}), (\ref{4.22}), and from 
the~inequality $\Im Z_{nk}(z)\ge \Im z,\,z\in\Bbb C^+$, that, 
for $z,\,k$ as above,
\begin{align}\label{4.23}
&|\theta_{nk}(Z_{nk}(z))|\le c(\tau)\varepsilon_{nk}^2Q^2(Z_{nk}(z))
\Big(|Z_{nk}(z)|+\varepsilon_{nk}Q(Z_{nk}(z))\Big)
\frac 1{|Z_{nk}(z)|^2}\notag\\
&\le c(\nu,\tau)\varepsilon_{nk}^2\frac 1{(\Im z)^4}Q^5(z)
\Big(1+\varepsilon_{nk}\frac 1{\Im z}Q(z)\Big).
\end{align}
Therefore (\ref{4.10}), (\ref{4.19}), (\ref{4.21}), and (\ref{4.23}) 
together yield the relation
\begin{equation}\label{4.24}
Z_{n1}(z)-z=-f_{n2}(Z_{n1}(z))+\dots
-f_{nk_n}(Z_{n1}(z))+r_n(z),\,\,\, z\in\Bbb C^+,
\end{equation}
where the~function 
$$
r_n(z):=\sum_{k=1}^{k_n}\big((F_{\widehat{\mu}_{nk}}(Z_{nk}(z))-Z_{nk}(z))
-(F_{\widehat{\mu}_{nk}}(Z_{n1}(z))-Z_{n1}(z))\big)+\sum_{k=1}^{k_n}
\theta_{nk}(Z_{nk}(z))
$$ 
is analytic in $\Bbb C^+$ and admits the~estimate
\begin{equation}\label{4.25}
|r_n(z)|\le  c(\nu,\tau)\varepsilon_n\frac 1{(\Im z)^4}Q^5(z)
\Big(1+\varepsilon_n\frac 1{\Im z}Q(z)\Big).
\end{equation}

From (\ref{4.25}) it is easy to see that
\begin{equation}\label{4.26}
|r_n(z)|\le c(\nu,\tau)\varepsilon_n^{1/20}
\end{equation}
in the~closed domain $D_n:=\{z\in\Bbb C^+:\varepsilon_n^{1/20}\le
\Im z\le \varepsilon_n^{-1/20},\,|\Re z|
\le \varepsilon_n^{-1/20}\}$.

We return to the~representation (\ref{4.14}) for the~functions $Z_{n1}(z)$. 
By (\ref{4.15}), (\ref{4.19a}), 
and the~vague compactness theorem (see \cite{Lo:1963}, p. 179), 
we conclude that there exists
a~subsequence $\{n'\}$ such that $Z_{n'1}(z+\gamma_{n'})\to Z(z)+az$ 
as $n'\to\infty$, uniformly on every compact set in $\Bbb C^+$, 
where $Z(z)\in\mathcal F$ and $a\ge 0$. 
Recalling the~assumption of the~theorem that $\nu_n\to \nu$ weakly and
(\ref{4.26}), we easily deduce from (\ref{4.24}) in the~limit $n'\to\infty$
that
\begin{equation}\label{4.27}
Z(z)+az-z=\int\limits_{\Bbb R}\frac{1+u(Z(z)+az)}{u-(Z(z)+az)}\,\nu(du),
\quad z\in\Bbb C^+.
\end{equation}
It is easy to see that $Z(iy)-iy=o(y)$ and the~integral on the~right-hand side
of (\ref{4.27}) is a~function which is $o(y)$ as $y\to\infty$ for $z=iy$. 
Therefore
we conclude that $a=0$. Thus the~relation (\ref{4.27}) holds with $a=0$.

Since $Z\in\Cal F$ has an~inverse $Z^{(-1)}$ defined on $\Gamma_{\alpha,\beta}$
with some positive $\alpha$ and $\beta$, it is easy to see that
the~equation (\ref{4.27}) has a~unique solution in the~set $\mathcal F$. 
Now suppose that $\{Z_{n1}(z+\gamma_n)\}_{n=1}
^{\infty}$ does not converge to $Z(z)$ on some compact set in $\Bbb C^+$. 
Then, as above there exists
a~subsequence $\{n''\}$ such that $Z_{n''1}(z+\gamma_{n''})\to Z^*(z)$ as 
$n''\to\infty$ on every compact set in $\Bbb C^+$, and $Z^*(z)\in\Cal F$, 
$Z^*(z)\not\equiv Z(z), \,z\in\Bbb C^+$. But $Z^*(z)$ is
a~solution of (\ref{4.27}). We arrive at a~contradiction. Hence  
$\{Z_{n1}(z+\gamma_n)\}_{n=1}^{\infty}$ converges to $Z(z)$ uniformly on 
every compact set 
in $\Bbb C^+$. The~relation (\ref{4.27}) implies that $Z(z)$ is infinitely 
divisible with parameters $(0,\nu)$, since we may rewrite (\ref{4.27})
via $z=Z^{(-1)}(w)$ for $w\in\Gamma_{\alpha,\beta}$ with some 
$\alpha,\beta>0$. Since 
$F_{\widehat{\mu}_{n1}}(Z_{n1}(z+\gamma_n))\to Z(z)$ uniformly 
on every compact set in $\Bbb C^+$, we see that $\widehat{\mu}_n\boxplus\delta_{-\gamma_n}$ 
converges weakly to a~p-measure $\widehat{\mu}$ such that 
$\phi_{\widehat{\mu}}=(0,\nu)$. Recalling the~definition of $a_n$, we
finally conclude that $\mu^{(n)}$ converges weakly 
to p-measure $\mu$ such that $\phi_{\mu}=(\alpha,\nu)$.

Hence the~sufficiency
of the~assumptions of Theorem~\ref{2.1,1}(b) is proved.

{\it Necessity}.
Denote $\mu_{nk}^s:=\mu_{nk}\boxplus\bar\mu_{nk},\,n\ge 1,\,k=1,\dots,k_n$.
By Proposition~\ref{3.5a}, we obtain the~convergence
\begin{equation}\label{4.28}
\mu^{(n,s)}:=\mu_{n1}^s\boxplus\mu_{n2}^s\boxplus\dots
\boxplus\mu_{nk_n}^s\to\mu^s
\quad\text{weakly as}\quad n\to\infty.
\end{equation}

For the~measures $\mu_{k,n}^s,\,n\ge 1,\,k=1,\dots,k_n$, relations (\ref{4.10}) 
and (\ref{4.11}) hold with the~functions 
$F_{\mu_{nk}^s}(z),\,n\ge 1,\,k=1,\dots,k_n$,
replacing $F_{\mu_{nk}}(z),\,n\ge 1,\,k=1,\dots,k_n$, and with some functions 
$Z_{nk,s}(z)\in\mathcal F,\,n\ge 1,\,k=1,\dots,k_n$, replacing  
$Z_{nk}(z),\,n\ge 1,\,k=1,\dots,k_n$. Rewrite (\ref{4.10}) in the~form
\begin{equation}\label{4.29}
F_{\mu_{n1}^s}(Z_{n1,s}(z))-z=F_{\mu_{n1}^s}(Z_{n1,s}(z))
-Z_{n1,s}(z)+\dots
+F_{\mu_{nk_n}^s}(Z_{nk_n,s}(z))-Z_{nk_n,s}(z) ,\quad z\in\Bbb C^+.
\end{equation}
By Proposition~\ref{3.7a}, the~measures $\mu_{nk}^s,\,k=1,\dots,k_n$, are
symmetric and $\mu^{(n,s)}:=\mu_{n1}^s\boxplus\mu_{n2}^s\boxplus\dots
\boxplus\mu_{nk_n}^s$ is symmetric as well. Since $F_{\mu_{nk}^s}(Z_{nk,s}(z))
=F_{\mu^{(n,s)}}(z),\,z\in\Bbb C^+$, and by Proposition~\ref{3.7a},
$F_{\mu^{(n,s)}}(iy),\,F_{\mu_{nk}^s}(iy),\,y>0$, assume imaginary values, 
we conclude that $Z_{nk,s}(iy),\,y>0,\,k=1,\dots,k_n$, assume imaginary 
values as well.
Hence, it is easy to see that $Z_{nk,s}(z),\,k=1,\dots,k_n$, 
admit representation (\ref{4.14}) with $d_{nk}=0,\,k=1,\dots,k_n$, 
and with finite nonnegative symmetric measures 
$\nu_{nk}^s,\,k=1,\dots,k_n$, respectively. 

Since
$$
F_{\mu_{nk}^s}(z)=z+\int\limits_{\Bbb R}
\frac{1+uz}{u-z}\,\sigma_{nk}^s(du),\quad z\in\Bbb C^+,\quad k=1,\dots,k_n,
$$
where $\sigma_{nk}^s$ is a~finite nonnegative measure, we deduce 
from (\ref{4.29}) that
\begin{align}\label{4.30}
\Im(F_{\mu_{n1}^s}(Z_{n1,s}(i))-i)
&=\sum_{k=1}^n\Im(F_{\mu_{nk}^s}(Z_{nk,s}(i))-Z_{nk,s}(i))\notag\\
&=\sum_{k=1}^n\Im Z_{nk,s}(i)\int\limits_{\Bbb R}
\frac{(1+u^2)\,\sigma_{nk}^s(du)}{u^2+(\Im Z_{nk,s}(i))^2}.
\end{align}

Let us show that the~measures $\mu_{n1}^s,\dots,\mu_{nk_n}^s$ are
infinitesimal. Indeed, we deduce from (\ref{4.9}) the~estimate
$$
-\Im \phi_{\mu_{nk}^s}(z)=-\Im \phi_{\widehat{\mu}_{nk}}(z)-
\Im \phi_{\overline{\widehat{\mu}}_{nk}}(z)\le c(\varepsilon_{nk}+
\overline{\varepsilon}_{nk})\le c\varepsilon_{nk},\quad k=1,\dots,k_n,
$$
for $|z-i|\le 1/4$, where $\overline{\varepsilon}_{nk}
:=\int_{\Bbb R}u^2/(1+u^2)\,\overline{\widehat{\mu}}_{nk}(du)
=\varepsilon_{nk}$. This implies
\begin{equation}\label{4.30a}
\int\limits_{\Bbb R}\frac {u^2}{1+u^2}\,\mu_{nk}^s(du)\Big/
\int\limits_{\Bbb R}\frac 1{1+u^2}\,\mu_{nk}^s(du)=
\Im(F_{\mu_{nk}^s}(i)-i)\le c\varepsilon_{nk},\quad k=1,\dots,k_n,
\end{equation}
as claimed.

The~bounds (\ref{4.30a}) and (\ref{4.8}) for the~functions 
$F_{\mu_{nk}^s}(z),\,n\ge n_0,\,k=1,\dots,k_n$, implies the~inequality
\begin{equation}\label{4.30b}
|Z_{nk,s}(i)|\le |Z_{nk,s}(i)-F_{\mu_{nk}^s}(Z_{nk,s}(i)|
+|F_{\mu^{(n,s)}}(i)|\le c\varepsilon_{nk}Q(Z_{nk,s}(i))
+|F_{\mu^{(n,s)}}(i)|
\end{equation}
for $n\ge n_0,\,k=1,\dots,k_n$. Since $Z_{nk,s}\in\Cal F$ and takes imaginary values for 
$z=iy,\,y>0$, we see that $|Z_{nk,s}(i)|=\Im Z_{nk,s}(i)\ge 1$. We note from this that
$Q(Z_{nk,s}(i))\le 2|Z_{nk,s}(i)|$ and,
by (\ref{4.28}), we easily conclude
from (\ref{4.30b}) that
$$
|Z_{nk,s}(i)|\le c(\mu^s),\quad n\ge n_0,\,\,k=1,\dots,k_n.
$$  
Moreover, 
$$
\Im(F_{\mu_{n1}^s}(Z_{n1,s}(i))-i)=\Im(F_{\mu^{(n,s)}}(i))-i)\to
\Im(F_{\mu^{(s)}}(i))-i),\quad n\to\infty. 
$$
Therefore we obtain from (\ref{4.30}) the~relation
\begin{equation}\label{4.31}
\sigma_{n1}^s(\Bbb R)+\dots+\sigma_{nk_n}^s(\Bbb R)\le c(\mu^s),
\quad n\to\infty. 
\end{equation}

Since $\mu_{nk}^s=\widehat{\mu}_{nk}\boxplus\bar{\widehat{\mu}}_{nk}$, 
we note, by definition of the free $\boxplus$-convolution (see Section~3),
that there exist functions $W_{nk}(z)\in\mathcal F$ such that
$F_{\mu_{nk}^s}(z)=F_{\widehat{\mu}_{nk}}(W_{nk}(z)),\,z\in\Bbb C^+$.
Therefore we have $\Im F_{\mu_{nk}^s}(i)-1=\Im F_{\widehat{\mu}_{nk}}
(W_{nk}(i))-1$. Rewrite this relation in the~form
\begin{align}\label{4.32}
\sigma_{nk}^s(\Bbb R)&= \Im W_{nk}(i)-1+\Im W_{nk}(i)\int\limits_{\Bbb R}
\ffrac{1+u^2}{(u-\Re W_{nk}(i))^2+(\Im W_{nk}(i))^2}\,\sigma_{nk}(du)\notag\\
&\ge \Im W_{nk}(i)\int\limits_{\Bbb R}  
\ffrac{1+u^2}{(u-\Re W_{nk}(i))^2+(\Im W_{nk}(i))^2}\,\sigma_{nk}(du).  
\end{align}
As in the~proof of (\ref{4.8}), we see that $F_{\mu_{nk}^s}(z)$ 
and $F_{\widehat{\mu}_{nk}}(z)$
tend to $z$ as $n\to\infty$ uniformly in $k=1,\dots,k_n$ and $|z-i|<1/2$. 
Hence $W_{nk}(i)\to i$ as $n\to\infty$ uniformly in $k=1,\dots,k_n$.  
Thus we obtain from (\ref{4.32}) that, for sufficiently large $n\ge n_0$,
\begin{equation}\label{4.33}
\sigma_{nk}(\Bbb R)\le 2\sigma_{nk}^s(\Bbb R),
\quad k=1,\dots,k_n.  
\end{equation}
Thus (\ref{4.31}) and (\ref{4.33}) imply
the~inequality
\begin{equation}\label{4.34}
\sigma_{n1}(\Bbb R)+\dots+\sigma_{nk_n}(\Bbb R)\le c(\mu^s),
\quad n\to\infty.
\end{equation}
By (\ref{4.2a}), (\ref{4.3a}) with $y=1$, and (\ref{4.6}), we note that 
$\sigma_{nk}(\Bbb R)\ge \varepsilon_{nk}/2,\,k=1,\dots,k_n$, 
for sufficiently large $n\ge n_0$ and we deduce from
(\ref{4.34}) the~upper bound
\begin{equation}\label{4.35}
\varepsilon_{n1}+\dots+\varepsilon_{nk_n}\le c(\mu^s),
\quad n\to\infty.
\end{equation}

Let us return to (\ref{4.10}) and (\ref{4.11}). 
Since $F_{\mu^{(n)}}(z)=F_{\widehat{\mu}_n}(z+a_n-b_n)$, where
$b_n=\sum_{k=1}^{k_n}a_{nk}$, we see that $F_{\mu^{(n)}}(z)=
F_{\widehat{\mu}_{nk}}(Z_{nk}(z+a_n-b_n)),\,z\in\Bbb C^+,\,k=1,\dots,k_n$.
Since $F_{\widehat{\mu}_{nk}}(z)$ tend to $z$ and 
$F_{\widehat{\mu}_{nk}}(Z_{nk}(z+a_n-b_n))$ tend to $F_{\mu}(z)$ 
as $n\to\infty$ 
uniformly in $k=1,\dots,k_n$ and on every compact set in $\Bbb C^+$,
we obtain that $\{Z_{nk}(z+a_n-b_n)\}_{n=1}^{\infty}$ converges
uniformly in $k=1,\dots,k_n$ and on every compact set in $\Bbb C^+$ 
to the~function $Z(z):=F_{\mu}(z)\in\mathcal F$. Using relations (\ref{4.10}) 
and (\ref{4.11}) with $z+a_n-b_n$ instead of $z$ and taking into account 
that the~measures $\mu_{n1},\dots,\mu_{nk_n}$
are infinitesimal and the~upper bound (\ref{4.35}) holds, we can 
repeate the~arguments which we used for the~proof of (\ref{4.24}). 
We arrive at the~following relation, for $z\in\Bbb C^+$, 
\begin{equation}\label{4.36}
Z_{n1}(z+a_n-b_n)-(z+a_n-b_n)=-f_{n2}(Z_1(z+a_n-b_n))-\dots
-f_{nk_n}(Z_1(z+a_n-b_n))+r_n(z),
\end{equation}
where $r_n(z)$ is analytic in $\Bbb C^+$ and $r_n(z)\to 0$ on every
compact set in $\Bbb C^+$.
As above, 
$\{Z_{nk}(z+a_n-b_n)\}_{n=1}^{\infty}$ converges uniformly in 
$k=1,\dots,k_n$ and on every 
compact set in $\Bbb C^+$ to $Z(z)\in\mathcal F$. 
Since, by (\ref{4.35}), 
the~sequence $\{\nu_n\}_{n=1}^{\infty}$ is tight in the~vague topology, 
there exists a~subsequence $\{n'\}$ such that there exists 
$\lim_{n'\to\infty}\nu_{n'}(\Bbb R)<\infty$ and $\{\nu_{n'}\}$ converges
to some finite nonnegative measure $\nu$ in the~vague topology. 
Now we conclude
from (\ref{4.36}) that $(a_{n'}-b_{n'}-\gamma_{n'})\to a'$ as $n'\to\infty$,
where $a'\in\Bbb R$, and the~following relation holds
\begin{equation}\label{4.37}
Z(z)=z+a'+b'Z(z)+\int\limits_{\Bbb R}\frac{1+uZ(z)}{u-Z(z)}\,\nu(du),
\quad z\in\Bbb C^+,
\end{equation}
with $b'=\lim_{n'\to\infty}\nu_{n'}(\Bbb R)-\nu(\Bbb R)\ge 0$. Recalling that 
$Z(z)\in\Cal F$, we easily conclude
from this relation that $b'=0$. Indeed, it is not difficult to see that
$Z(iy)-iy=o(y)$ and the~integral in (\ref{4.37}) for $z=iy$ is $o(y)$ as $y\to+\infty$.
Comparing a~behaviour of all members in (\ref{4.37}), we obtain 
the~desired result.

We shall show that $\{\nu_{n}\}$ 
converges to the~measure $\nu$ in the~vague topology. Assume to the~contrary
that there exists a~subsequence $\{n''\}$ such that there exists 
$\lim_{n''\to\infty}\nu_{n''}(\Bbb R)<\infty$ and $\{\nu_{n''}\}$ 
converges in the~vague topology to some finite measure $\nu''\not\equiv \nu$.  
Then $(a_{n''}-b_{n''}-\gamma_{n''})\to a''$ as $n''\to\infty$,
and (\ref{4.37}) holds with $a''$ replacing $a'$ and $\nu''$ replacing
$\nu$. Comparing relations (\ref{4.37}), we deduce the~relation
$$
a'+\int\limits_{\Bbb R}\frac{1+uz}{u-z}\,\nu(du)=
a''+\int\limits_{\Bbb R}\frac{1+uz}{u-z}\,\nu''(du),
\quad z\in\Bbb C^+.
$$
Applying Stieltjes-Perron inversion formula (see Section~3), we get
that $\nu=\nu''$ and then $a'=a''$, a~contradiction. Since, as above,
$\lim_{n\to\infty}\nu_n(\Bbb R)=\nu(\Bbb R)$, we finally conclude 
that $\{\nu_n\}$ converges
to the~measure $\nu$ weakly. In addition, $a_n-b_n-\gamma_n$ tends to
some real constant as $n\to\infty$. It remains to note that, by
the~relation $Z(z)=F_{\mu}(z), \,z\in\Bbb C^+$, we see from (\ref{4.37})
that the~limit measure $\mu$ is infinitely divisible with parameters
$(a',\nu)$.

This proves the~necessity of the~assumptions of Theorem~2.1(b) and
thus Theorem~2.1.
$\square$

{\bf Proof of Corollary~\ref{2.2,1a}.}
In order to prove Corollary~\ref{2.2,1a} using Theorem~\ref{2.2,1} 
we need to show that if $\mu_n^{k_n*}$ converges weakly to some 
$\mu^*\in\Cal M$ or $\mu_n^{k_n\boxplus}$ converges 
to some $\mu^{\boxplus}\in\Cal M$, then $\mu_n$ are infinitesimal.  
The~first assertion is a~well-known fact (see \cite{Lo:1963}).
It remains to prove the~second assertion only. 
By Proposition~\ref{3.4a},
$k_n\phi_{\mu_n}(z)$
converges uniformly on compact subsets of $\Gamma_{\alpha,\beta}$, 
with some $\alpha,\beta>0$, to the~function $\phi_{\mu^{\boxplus}}(z)$
and $k_n\phi_{\mu_n}(iy)=o(y)$ 
uniformly in $n$ as $y\to+\infty$. 
Hence $\phi_{\mu_n}(z)\to 0$ uniformly on compact subsets of 
$\Gamma_{\alpha,\beta}$, as $n\to\infty$, and $\phi_{\mu_n}(iy)=o(y)$ 
uniformly in $n$ as $y\to+\infty$. By Proposition~\ref{3.4a}, 
$\mu_n$ converges weakly to $\delta_0$ as $n\to\infty$. 
Therefore the~p-measures $\mu_n$ are infinitesimal.
$\square$

\section{The~class $\mathcal L_{\boxplus}$ of infinitely divisible 
limits for free sums}

In this section we shall prove Theorem~\ref{2.7,1} and Theorem~\ref{2.10,1}.

Let $\{\mu_n\}_{n=1}^{\infty}$ be a~sequence of measures in $\mathcal M$
and let $\{a_n\}_{n=1}^{\infty}$ and $\{b_n\}_{n=1}^{\infty}, b_n>0$, be
sequences of real numbers. In Section~5 we denote by $\mu_{nk}:n\ge 1,
1\le k\le n$, the~p-measures such that $\mu_{nk}(S):=\mu_k(b_nS)$
for every Borel set $S\in\Bbb R$. Consider the~sequence of p-measures
$\{\mu^{(n)}:=\delta_{-a_n}\boxplus\mu_{n1}\boxplus\dots\boxplus
\mu_{nn}\}_{n=1}^{\infty}$. Recall that we denote by 
$\widehat{\mu}_{nk},n\ge 1,k=1,\dots,n$, the~p-measures 
defined by (\ref{2.9}). By $\nu_n,\,n=1,\dots$, we denote
the~measures defined in (\ref{2.10}).

In the~first step we shall prove Theorem~\ref{2.7,1}. First we prove
the~following auxiliary proposition.
\begin{proposition}\label{proposition 5.1}
If $\mu^{(n)}\to\mu$ weakly as $n\to\infty$, where $\mu$ is 
a~nondegenerate p-measure, and if 
the~measures $\mu_{n1},\dots \mu_{nn}$ are infinitesimal, then 
$b_n\to\infty$ and $b_{n+1}/b_n\to 1$ as $n\to\infty$ hold.
\end{proposition}
{\bf Proof.}
Assume to the~contrary that there exists an~infinite sequence $\{n'\}$ 
such that the~numbers $b_{n'}$ remain bounded. Without loss of generality,
we may assume that the~numbers $b_{n'}$ converge to some finite $b>0$ as 
$n'\to\infty$. By Proposition~\ref{3.4a},
and since $\mu_{nk},\,n\ge 1,\,k=1,\dots,n$, are infinitesimal, for some 
$\alpha,\beta>0$ and $z\in\Gamma_{\alpha,\beta}$, $\phi_{\mu^{(n)}}(z)$
converges, that is,
\begin{equation}\label{5.1}
\phi_{\mu^{(n)}}(z)=-a_n+\phi_{\mu_{n1}}(z)+\dots+\phi_{\mu_{nn}}(z)=
-a_n+\frac 1{b_n}\phi_{\mu_1}(b_nz)+\dots+\frac 1{b_n}
\phi_{\mu_n}(b_nz)\to \phi_{\mu}(z) 
\end{equation} 
as $n\to\infty$. By Proposition~\ref{3.4a},
and since $\mu_{nk},\,n\ge 1,\,k=1,\dots,n$, are infinitesimal,
we conclude for every $k$, $\alpha,\beta>0$, and $z\in\Gamma_{\alpha,\beta}$, 
that $\frac 1{b_{n'}}\phi_{\mu_k}(b_{n'} z)\to 0$ as $n'\to\infty$. 
Hence $\phi_{\mu_k}(z)\equiv 0,\,z\in\Gamma_{\alpha,\beta},\,k=1,\dots$, 
and $\phi_{\mu^{(n)}}(z)
=-a_n$ for $z\in \Gamma_{\alpha,\beta}$. Since $\mu$ is nondegenerate,
we arrive at a~contradiction, using Proposition~\ref{3.4a}.

Since $\mu_{nk},\,n\ge 1,\,k=1,\dots,n$, are infinitesimal, we note that
$
\mu^{(n+1,n)}:=\delta_{-a_{n+1}}\boxplus\mu_{n+1,1}\boxplus\dots\boxplus
\mu_{n+1,n}
$
also converge to $\mu$ weakly as $n\to\infty$. Denote the~distribution 
function of $\mu^{(n)}$ by $F_n(x)$. Then 
the~distribution function of $\mu^{(n+1,n)}$ has 
the~form $F_n(b'_nx+a'_n)$, where $b'_n=b_{n+1}/b_n$ and 
$a'_n=(b_{n+1}/b_n)a_{n+1}-a_n$.
According to Proposition~\ref{3.8a}, this implies that $(b_{n+1}/b_n)\to 1$
as $n\to\infty$, thus proving the~proposition.
$\square$

{\bf Proof of Theorem~\ref{2.7,1}.}
The~``if'' assertion follows from Theorem~\ref{2.1,1} by choosing the~measures
$\mu^{(n)}=\delta_{-a_n}\boxplus\mu_{n1}\boxplus\dots\boxplus\mu_{nn}$,
where $\mu_{nk},\,k=1,\dots,n$, are defined at the~beginning
of this section with $b_n=b_n'$.

Now we assume that the~distributions $\mu^{(n)}$ 
converge to a~limit distribution $\mu$ 
for some choice of the~sequence of constants $b_n>0$ and $a_n$, 
and assume that the~measures $\mu_{nk},\,n\ge 1,\,k=1,\dots,n$, are infinitesimal. 
By Theorem~\ref{2.1,1}, $\mu$ is 
an~infinitely divisible probability measure and its Voiculescu's 
transform $\phi_{\mu}(z)$ has the~form (\ref{2.6}).
If we prove that $(b_n/b'_n)\to 1$ as $n\to\infty$,
then, by Proposition~\ref{3.8a}, the~sequence of $\mu^{(n)}$ (where
$\mu_{nk},\,k=1,\dots,n$, defined using $b_n=b_n'$)
converges to $\mu$ weakly too. Hence (\ref{2.12}) would follow from 
Theorem~\ref{2.1,1} and the ''only if'' assertion would hold.

Consider the~p-measures $\mu^{(n)s}:=\mu_{n1}^s\boxplus\dots
\boxplus\mu_{nn}^s$. Since the~p-measure $\mu_{n1},\dots,\mu_{nn}$ 
are infinitesimal, as we proved in Section~4, the~p-measures 
$\mu_{n1}^s,\dots,\mu_{nn}^s$ are infinitesimal as well.
By Proposition~\ref{3.5a}, 
$\mu^{(n)s}\to \mu^s:=\mu\boxplus \bar{\mu}$ weakly as $n\to\infty$. 
This measure is infinitely divisible
and its Voiculescu's transform $\phi_{ \mu^s}(z)$ admits 
a~representation (\ref{2.6}) with parameters $(0,\nu_s)$, where $\nu_s:=
\nu+\bar\nu$. Moreover, by Theorem~\ref{2.1,1}, $\nu_{ns}
\to\nu_s$ weakly, as $n\to\infty$, where the~$\nu_{ns}$ are defined by
$$
\nu_{ns}(S):=\sum_{k=1}^n\int\limits_{b_n S}\frac{u^2}{b_n^2+u^2}\,
\mu_k^s(du)
$$ 
for every Borel set $S$. Here $b_nS:=\{b_nx:x\in S\}$. 
Using the~relation $\nu_s(\Bbb R)=2\nu(\Bbb R)$
and (\ref{2.10a}), we obtain
\begin{equation}\label{5.2}
\Delta_n=\sum_{k=1}^n\Big(\int\limits_{\Bbb R}\frac{u^2}{b_n^2+u^2}
\,\mu_k^s(du)-\int\limits_{\Bbb R}\frac{u^2}{(b'_n)^2+u^2}\,\mu_k^s(du)
\Big)\to 0.
\end{equation}
By the~assumption, $\mu$ is a~nondegenerate p-measure. 
Therefore, there exists an~$a>0$ such that $2\delta=\nu_s((-a,a))>0$
and, hence, for sufficiently large $n\ge n_a$,
$$
\sum_{k=1}^n\int\limits_{-ab_n}^{ab_n}\frac{u^2}{b_n^2+u^2}\,\mu_k^s(du)
>\delta>0.
$$   
Hence,
\begin{align}
|\Delta_n|&=|b_n^2-(b'_n)^2|\sum_{k=1}^n\int\limits_{\Bbb R}
\frac{u^2}{(b_n^2+u^2)((b'_n)^2+u^2)}\,\mu_k^s(du)\notag\\
&\ge\frac{|b_n^2-(b'_n)^2|}{(b'_n)^2+a^2b_n^2}
\sum_{k=1}^n\int\limits_{-ab_n}^{ab_n}\frac{u^2}{b_n^2+u^2}
\,\mu_k^s(du)\ge \ffrac{|(b_n/b'_n)^2-1|}{1+a^2(b_n/b'_n)^2}\delta.\notag
\end{align}
By (\ref{5.2}), this yields $b_n/b'_n\to 1$ and the~proof
is complete.
$\square$

Now we shall prove Theorem~\ref{2.10,1}.

{\bf Proof of Theorem~\ref{2.10,1}.} {\it Sufficiency}. 
First we show that $\phi_{\mu}(z)$ admits an~analytic continuation
on $\Bbb C^+$. Indeed, the~function $\phi_{\mu}(z)$ is regular on some domain  
$\Gamma_{\alpha,\beta}$ (see Section~2). Let us assume that there are 
singular points  
of $\phi_{\mu}(z)$ on the~boundary of this domain. Let $z_0$ be one of such 
points with the~largest modulus. By the~definition of $\phi_{\mu}(z)$,
it is easy to see that $|z_0|<\infty$.
By the~assumption and by Voiculescu's relation~(\ref{2.4}), 
we have
\begin{equation}\label{5.3}
\phi_{\mu}(z)=\gamma\phi_{\mu}(z/\gamma)
+\phi_{\gamma}(z),\quad z\in\Gamma_{\alpha,\beta},
\end{equation}
for $0<\gamma<1$. Here $\phi_{\gamma}(z):=\phi_{\mu_{\gamma}}(z)$.
Hence
$$
\phi_{\gamma}(z)=(1-\gamma)\phi_{\mu}(z)+\gamma(\phi_{\mu}(z)
-\phi_{\mu}(z/\gamma))=
(1-\gamma)\phi_{\mu}(z)+\gamma\int\limits_{z/\gamma}^z\phi'_{\mu}(\zeta)\,d\zeta,
\quad z\in\Gamma_{\alpha,\beta}.
$$
By Proposition~\ref{3.6a} and the~property $|\phi_{\mu}(z)|=o(|z|)$ as 
$z\to\infty,\,z\in\Gamma_{\alpha,\beta}$, we have the~relation 
$\phi_{\gamma}(z)\to 0$
for $z\in\Gamma_{\alpha,\beta}$ and $\phi_{\gamma}(iy)=o(y),\, y\to\infty$, 
uniformly in $\gamma\to 1$.
If $\gamma\to 1$, then by Proposition~\ref{3.4a}, $\mu_{\gamma}\to\delta_0$ 
weakly and, 
by Proposition~\ref{3.6aa}, $\phi_{\mu_{\gamma}}(z)$ is regular in the~domain 
$\Gamma_{2\alpha,\beta/2}$ for $\gamma$ close to $1$.  
The~functions $\phi_{\mu}(z/\gamma)$ and $\phi_{\gamma}(z)$
are regular on $\overline {\Gamma}_{\alpha,\Im z_0}$, therefore $\phi_{\mu}(z)$ 
is regular at the~point $z_0$, a~contradiction. Hence our assertion holds. 
Note that the~function $\phi_{\gamma}(z)$ admits an~analytic continuation
on $\Bbb C^+$ for every $\gamma\in(0,1)$ as well and the~relation (\ref{5.3})
holds for all $\gamma\in(0,1)$ and $z\in\Bbb C^+$ for such functions.
We again denote these functions, defined on $\Bbb C^+$, by $\phi_{\mu}(z)$
and $\phi_{\gamma}(z)$.

Consider the~p-measures $\mu_k,\,k=1,\dots$, determined via
$$
\phi_{\mu_k}(z):=\frac 1{\pi_k}\phi_{\gamma_k}(\pi_k z)
=\frac 1{\pi_k}\phi_{\mu}(\pi_k z)
-\frac 1{\pi_{k-1}}\phi_{\mu}(\pi_{k-1}z),\quad z\in\Bbb C^+,
$$
where $\pi_k=\prod_{l=1}^k\gamma_l$ and $\gamma_l=1-1/(l+1),\,\pi_0:=1$.

Voiculescu's transform of the~p-measure $\mu^{(n)}:=\mu_{n1}
\boxplus\dots\boxplus\mu_{nn}$ with $\mu_{nk},\,k=1,\dots,n$, 
defined at the~beginning of this section using $b_n:=1/\pi_n$,
has the~form
\begin{equation}\label{5.4}
\phi_{\mu^{(n)}}(z)=\sum_{k=1}^n\Big(\frac{\pi_n}{\pi_k}
\phi_{\mu}\Big(\frac{\pi_k}{\pi_n} z\Big)-\frac{\pi_n}
{\pi_{k-1}}\phi_{\mu}\Big(\ffrac{\pi_{k-1}}{\pi_n}z\Big)\Big)
=\phi_{\mu}(z)-\pi_n\phi_{\mu}\Big(\frac z{\pi_n}\Big),
\quad z\in\Bbb C^+.
\end{equation}
From (\ref{5.4}) and Proposition~\ref{3.4a} it follows that $\mu^{(n)}\to\mu$
weakly as $n\to\infty$. It is not difficult to verify that 
the~functions $\phi_{\mu_{nk}}(z)$ converge to zero 
on every compact set of $\Bbb C^+$ uniformly in $k,\,1\le k\le n$. 
This implies that $|F_{\mu_{nk}}(i)-i|\to 0$ as $n\to\infty$
uniformly in $k=1,\dots,n$, and hence $|G_{\mu_{nk}}(i)+i|\to 0$
as as $n\to\infty$ uniformly in $k=1,\dots,n$. Finally note
that the~p-measures $\mu_{n1},\dots,\mu_{nn}$ are
infinitesimal. This follows directly from formula (\ref{4.3a}) 
for the~measures $\mu_{nk},\,k=1,\dots,n$.

{\it Necessity}. 
Let $\mu\in\mathcal L_{\boxplus}$. This means that there exists
a~sequence of p-measures $\{\mu_n\}$ such that 
for some suitably chosen sequences of constants  $\{a_n\}$ and 
$\{b_n\},\,b_n>0$, the~sequence of the~measures 
$\{\mu^{(n)}=\delta_{-a_n}\boxplus\mu_{n1}\boxplus\dots\boxplus\mu_{nn}\}$
converges weakly to a~limit measure $\mu$ and that the~measures
$\mu_{nk},\,n\ge 1,\,k=1,\dots,n$, are infinitesimal. Hence condition (\ref{4.1})
holds for the~measures $\mu_{nk},\,k=1,\dots,n$. 
By Proposition~\ref{3.6aa}, for every $\alpha,\beta>0$ Voiculescu's
transforms $\phi_{\mu_{nk}}(z),\,k=1,\dots,n$, are defined on
$\Gamma_{\alpha,\beta}$ and $\phi_{\mu_{nk}}(\Gamma_{\alpha,\beta})
\subset \Bbb C^-\cup\Bbb R,\,k=1,\dots,n$. Expressing our assumptions
in terms of Voiculescu's transforms, we obtain, by Proposition~\ref{3.4a},
that 
\begin{equation}\label{5.5}
\phi_{\mu^{(n)}}(z)=-a_n+\frac 1{b_n}\phi_{\mu_1}(b_nz)+\dots
+\frac 1{b_n}\phi_{\mu_n}(b_nz)\to \phi_{\mu}(z),\quad n\to\infty,
\end{equation}
uniformly on compact subsets of $\Gamma_{\alpha,\beta}$ and
$\phi_{\mu^{(n)}}(iy)=o(y)$ uniformly in $n$ as $y\to+\infty$.
The~function $\phi_{\mu}(z)$ in (\ref{5.5}) is regular in $\Bbb C^+$
(being Voiculescu's transform of an~infinitely divisible p-measure).
According to Proposition~5.1, for every given $\gamma\in(0,1)$
there exists an~$m=m(n), m<n$, such that $b_m/b_n\to\gamma$
as $n\to\infty$. Rewrite $\phi_{\mu^{(n)}}(z)$ in the~form
\begin{align}\label{5.6}
\phi_{\mu^{(n)}}(z)=\Big(-\frac{b_m}{b_n}a_m&+\frac 1{b_n}\phi_{\mu_1}(b_nz)
+\dots+\frac 1{b_n}\phi_{\mu_m}(b_nz)\Big)\notag\\
&+\Big(\frac{b_m}{b_n}a_m-a_n+\frac 1{b_n}\phi_{\mu_{m+1}}(b_nz)+\dots
+\frac 1{b_n}\phi_{\mu_n}(b_nz)\Big),
\quad z\in\Gamma_{\alpha,\beta}.
\end{align}
Since, by (\ref{5.5}), 
\begin{equation}\notag
-a_m+\frac 1{b_m}\phi_{\mu_1}(b_mz)+\dots+\frac 1{b_m}
\phi_{\mu_m}(b_mz)\to \phi_{\mu}(z),\quad m\to\infty,
\end{equation}
uniformly on compact subsets of $\Gamma_{\alpha,\beta}$ and
$\phi_{\mu^{(m)}}(iy)=o(y)$ uniformly in $n$ as $y\to+\infty$, we have 
\begin{equation}\label{5.7}
\phi_{m,n}(z):=-\frac{b_m}{b_n}a_m+\frac 1{b_n}\phi_{\mu_1}(b_nz)
+\dots+\frac 1{b_n}\phi_{\mu_m}(b_nz)\to \gamma\phi_{\mu}
\big(\frac z{\gamma}\big),\quad n\to\infty,
\end{equation}
uniformly on compact subsets of $\Gamma_{\alpha,\beta}$ and
$\phi_{m,n}(iy)=o(y)$ uniformly in $n$ as $y\to+\infty$.
Therefore the~second bracket in (\ref{5.6}) converges to the~regular
function $\phi_{\gamma}(z):=\phi_{\mu}(z)-\gamma\phi_{\mu}(z/\gamma)$ 
on $\Bbb C^+$. By Proposition~\ref{3.4a}, $\phi_{\gamma}(z)$
is Voiculescu's transform of some measure $\mu_{\gamma}\in\Cal M$.
Thus we conclude that for every $\gamma\in(0,1)$ $\phi_{\mu}(z)=
\gamma\phi(z/\gamma)+\phi_{\gamma}(z)$ for all $z\in\Bbb C^+$.

Hence, the~theorem is proved.
$\square$

\section{ Estimates of convergence in the Free Central 
Limit Theorem}

In this section we prove Theorem~\ref{2.3,1}, Proposition~\ref{2.4,1}, 
Theorem~\ref{2.5,1}, and Theorem~\ref{2.6,1}.
In the sequel we denote $c_1,c_2\dots$ explicit positive absolute
constants.

{\bf Proof of Theorem~\ref{2.3,1}.} 
Denote $\mu^{(n)}:=\mu_n^{n\boxplus}$. By Corollary~\ref{3.3b}, 
$G_{\mu^{(n)}}(z)=1/F_{\mu^{(n)}}(z)$, $z\in\Bbb C^+$, 
where $F_{\mu^{(n)}}(z):=F_{\mu}(Z(\sqrt nz))/\sqrt n$.
In this formula $Z(z)\in\mathcal F$ is the~solution of 
equation~(\ref{3.8}). Consider the~functions 
$S(z):=\frac 12(z+\sqrt{z^2-4})$
and $S_n(z):=Z(\sqrt nz)/\sqrt n$ for $z\in\Bbb C^+$.  
Note that $1/S(z)=G_{\mu_w}(z)$, where $w$ denotes Wigner semicircle 
measure. Since $S_n\in\Cal F$, we see by Remark~\ref{3.0} that there exists
a p-measure $\nu^{(n)}$ such that $1/S_n(z)=G_{\nu^{(n)}}(z)$. 

We obtain the~estimate (\ref{2.9,a}) for $n\ge n_2$, 
where $n_2:=[c_1(|m_3(\mu)|^2+m_4(\mu))]$ with a~sufficiently large 
positive absolute constant $c_1$. For $n\le n_2$ (\ref{2.9,a}) holds 
obviously. Using (\ref{2.1}), we may write
\begin{align}\label{6.4,1}
Z(z)G_{\mu}(Z(z))&=1+\frac 1{Z^2(z)}+\frac 1{Z^2(z)}
\int\limits_{\Bbb R}\frac{u^3\,\mu(du)}{Z(z)-u}\notag\\&=
1+\frac 1{Z^2(z)}+\frac {m_3(\mu)}{Z^3(z)}+
\frac 1{Z^3(z)}\int\limits_{\Bbb R}\frac{u^4\,\mu(du)}{Z(z)-u},
\quad z\in\Bbb C^+.
\end{align}
Equation (\ref{3.8}) may be rewritten as
\begin{equation}\label{6.5,1}
G_{\mu}(Z(z))\Big(Z(z)-z\Big)=(n-1)(1-Z(z)G_{\mu}(Z(z))),
\quad z\in\Bbb C^+.
\end{equation}
By (\ref{6.4,1}) and the~definition of $S_n(z)$, 
(\ref{6.5,1}) may be reformulated as
\begin{align}\label{6.6,1}
\Big(1&+\frac 1{Z^2(\sqrt n z)}+\frac {m_3(\mu)}{Z^3(\sqrt nz)}+
\frac 1{Z^3(\sqrt nz)}\int\limits_{\Bbb R}\frac{u^4\,\mu(du)}
{Z(\sqrt nz)-u}\Big)(S_n(z)-z)\notag\\
&=-\frac{n-1}n\Big(\frac 1{S_n(z)}+\frac {m_3(\mu)}{S_n(z)Z(\sqrt nz)}
+\frac 1{S_n(z)Z(\sqrt nz)}
\int\limits_{\Bbb R}\frac{u^4\,\mu(du)}{Z(\sqrt nz)-u}\Big)
\end{align}
for $z\in\Bbb C^+$. Rewrite (\ref{6.6,1}) in the~form
\begin{equation}\label{6.7,1}
(1+r_{n1}(z))(S_n(z)-z)=-\Big(1-\frac 1n\Big)\frac 1{S_n(z)}(1+r_{n2}(z)),
\end{equation}
where $r_{n1}(z)$ and $r_{n2}(z)$ are analytic functions on $\Bbb C^+$ 
which, by the~inequality $\Im Z(\sqrt nz)\ge\sqrt n\Im z,\,
z\in\Bbb C^+$, (compare with (\ref{3.2,a})), admit the~estimates
\begin{equation}\label{6.8,1}
|r_{n1}(z)|\le\frac 1{(\Im z\sqrt n)^2}+\frac{|m_3(\mu)|}
{(\Im z\sqrt n)^3}+\frac{m_4(\mu)}{(\Im z\sqrt n)^4},
\quad 
|r_{n2}(z)|\le\frac{|m_3(\mu)|}{\Im z\sqrt n}+
\frac{m_4(\mu)}{(\Im z\sqrt n)^2},\quad z\in\Bbb C^+. 
\end{equation}

Introduce for every $\alpha>0$, $\Bbb C^+_{\alpha}:=\{z\in\Bbb C:\,
\Im z> \alpha\}$ and $D_{\alpha}:=\{z\in\Bbb C: \alpha\le\Im z\le 1,
\,|\Re z|\le 4\}$.

By (\ref{6.8,1}), $|r_{n1}(z)|+|r_{n2}(z)|\le 1/10$ for $z\in\Bbb C^+_{a/2}$, 
where $a=:c_2(|m_3(\mu)|+m_4^{1/2}(\mu))/\sqrt n$ and $c_2>0$ 
is a~sufficiently large absolute constant. Therefore
we conclude from (\ref{6.7,1}) that 
\begin{equation}\label{6.8a,1}
10^{-1}\le |S_n(z)|\le 10,\quad z\in D_a.
\end{equation}


From (\ref{6.7,1}) we see that 
the~function $S_n(z)$ satisfies the~approximate functional equation
\begin{equation}\label{6.11,1}
S_n(z)-z=-\frac1{S_n(z)}+\frac{r_{n3}(z)}{S_n(z)} ,
\end{equation}
for $z\in\Bbb C_{a/2}^+$, where
$$
r_{n3}(z):=1-\Big(1-\frac 1n\Big)\frac{1+r_{n2}(z)}{1+r_{n1}(z)}.
$$
Here $r_{n3}(z)$ is an~analytic function on $z\in\Bbb C_{a/2}^+$ which
is bounded as follows
\begin{equation}\label{6.11a,1}
|r_{n3}(z)|\le 2\Big(\frac 1n+|r_{n1}(z)|+|r_{n2}(z)|\Big),
\quad z\in\Bbb C_{a/2}^+.
\end{equation}
Recalling the~definition of the~functions $r_{n1}(z)$ and $r_{n2}(z)$,
we obtain with the~help of (\ref{6.8a,1}) and the~inequality $|S_n(z)|\ge 1,
\,z\in\overline{\Bbb C}^+_1$,
\begin{equation}\label{6.11b,1}
|r_{n1}(z)|\le \frac 1{(\sqrt n|S_n(z)|)^2}+\frac{|m_3(\mu)|}
{(\sqrt n|S_n(z)|)^3}+\frac{m_4(\mu)}{(\sqrt n|S_n(z)|)^3\sqrt n \Im z}
\le 10^3\Big( \frac 1{n}+\frac{|m_3(\mu)|}{n^{3/2}}+\frac{m_4(\mu)}
{n^2\Im z}\Big)
\end{equation}
and
$$
|r_{n2}(z)|\le \frac{|m_3(\mu)|}{\sqrt n|S_n(z)|}
+\frac{m_4(\mu)}{n|S_n(z)|\Im z}\le 10\Big(\frac{|m_3(\mu)|}{\sqrt n}+
\frac{m_4(\mu)}{n\Im z}\Big)
$$
for $z\in D_a\cup \overline{\Bbb C}_1^+$. 
Applying these estimates to (\ref{6.11a,1}), we finally have
\begin{equation}\label{6.12,1}
|r_{n3}(z)|\le 
3\cdot 10^3\Big(\frac{|m_3(\mu)|}{\sqrt n}+\frac{m_4(\mu)} 
{n\Im z}+\frac 1n\Big),
\quad z\in D_a\cup\overline{\Bbb C}_1^+. 
\end{equation}
Solving equation (\ref{6.11,1}), we see that
$$
S_n(z)=\frac 12\Big(z\pm\sqrt{\rho_n(z)}\Big),\quad
z\in\Bbb C^+_{a/2},
$$ 
where $\rho_n(z):=z^2-4+4r_{n3}(z)$. Note that the~function $\rho_n(z)$
is non-zero on the~half-plane $\Bbb C^+_{a/2}$. Indeed, let $\rho_n(w)=0$ 
for some $w\in\Bbb C^+_{a/2}$. Then, by (\ref{6.11,1}), 
$S_n^2(w)-wS_n(w)=-w^2/4$ and
we have $S_n(w)=w/2$. But the~function $S_n(z)$ satisfies the~inequality
$\Im S_n(z)\ge \Im z,\,z\in\Bbb C^+$, a~contradiction. We define 
the~function $\sqrt{\rho_n(z)}$ on $\Bbb C^+_{a/2}$, taking the~branch
of $\sqrt{\rho_n(z)}$ such that $\sqrt{\rho_n(i)}\in\Bbb C^+$.
Since $S_n(z)\in\mathcal N$, we see that $S_n(z)=\frac 12
\Big(z+\sqrt{\rho_n(z)}\Big)$ for $z\in\Bbb C^+_{a/2}$.

For $z\in\Bbb C^+_{a/2}$, using the~previous formula for $S_n(z)$ and 
$S(z)=\frac 12(z+\sqrt{z^2-4})$, we write
\begin{equation}\label{6.13,1}
\frac 1{S_n(z)}-\frac 1{S(z)}=\frac{S(z)-S_n(z)}{S(z)S_n(z)}=
\frac 1{S(z)S_n(z)}\cdot\frac{2r_{n3}(z)}{\sqrt{z^2-4}
+\sqrt{z^2-4+4r_{n3}(z)}}.
\end{equation}
Since, for $z\in \Bbb C, \,0<\Im z\le 1$, $|z^2-4|\ge m(z):=\max\{\Im z,
((\Re z)^2-5)_+\}$, where for $x\in\Bbb R, (x)_+:=\max\{0,x\}$, 
we obtain from (\ref{6.12,1}) the~following inequality
\begin{equation}\label{6.14,1}
\Big|\frac{r_{n3}(z)}{z^2-4}\Big|\le \frac {3\cdot 10^3}{m(z)}
\Big(\frac {|m_3(\mu)|}{\sqrt n}
+\frac {2m_4(\mu)}{n\Im z}\Big)\le\frac 1{10},\quad z\in D_a\cup
\{z\in\Bbb C:\Im z=1\}.
\end{equation}
Hence we get, for $z\in D_a$ or for $\Im z=1$,
$$
|\sqrt{z^2-4}+\sqrt{z^2-4+4r_{n3}(z)}|=\sqrt{|z^2-4|}\Big|1
+\sqrt{1+4r_{n3}(z)/(z^2-4)}\Big|\ge \sqrt{|z^2-4|}.
$$
Using this estimate we deduce from (\ref{6.12,1}) and (\ref{6.13,1}),  
for $z\in(D_a\cup\{z\in\Bbb C:\Im z=1\})$, 
\begin{equation}\label{6.16,1}
\Big|\frac 1{S_n(z)}-\frac 1{S(z)}\Big|\le 2\frac{|r_{n3}(z)|}
{|\sqrt{z^2-4}|}\frac 1{|S(z)||S_n(z)|}
\le \frac {6\cdot 10^3}{\sqrt{m(z)}}\Big(\frac {|m_3(\mu)|}{\sqrt n}
+\frac {2m_4(\mu)}{n\Im z}\Big)\frac 1{|S(z)||S_n(z)|}.
\end{equation}

Recall that $1/S(z)=G_{\mu_w}(z)$ and 
$1/S_n(z)=G_{\nu^{(n)}}(z)$, where
$\nu^{(n)}$ is a~p-measure. 

Since, for $u\in\Bbb R$, $m(u+i)=\max\{1,(u^2-5)_+\}$, $|S_n(u+i)|\ge 1$,
and $|S(u+i)|\ge \frac 12 \sqrt{1+((u-4)_+)^2}$, we conclude, 
using (\ref{6.16,1}),
\begin{align}\label{6.17,1}
\int\limits_{\Bbb R}|G_{\mu_w}(u+i)&-G_{\nu^{(n)}}(u+i)|\,du
\le c\Big(\frac {|m_3(\mu)|}{\sqrt n}+\frac {m_4(\mu)}n\Big)
\int\limits_{\Bbb R}\frac{du}{1+u^2}\notag\\
&\le\frac c{\sqrt n}\Big(|m_3(\mu)|+\frac {m_4(\mu)}{\sqrt n}\Big)
\le\frac c{\sqrt n}\Big(|m_3(\mu)|+(m_4(\mu))^{1/2}\Big),\quad n\ge n_0.
\end{align} 

Since, for $z\in D_a$, $\sqrt{m(z)}\ge \sqrt{\Im z}$, $|S_n(z)|\ge 1/10$,  
and $|S(z)|\ge 1/10$, we obtain from (\ref{6.16,1}), for $x\in[-2,2]$,
\begin{align}\label{6.18,1} 
\int\limits_a^1|G_{\mu_w}(x+iu)&-G_{\nu^{(n)}}(x+iu)|\,du\le
c\int\limits_a^1\Big(\frac {|m_3(\mu)|}{\sqrt{nu}}
+\frac {m_4(\mu)}{nu^{3/2}}\Big)\,du\notag\\
&\le \frac c{\sqrt n}\Big(|m_3(\mu)|+\frac {m_4(\mu)}{\sqrt {na}}\Big)
\le \frac c{\sqrt n}\Big(|m_3(\mu)|+(m_4(\mu))^{1/2}\Big),\quad n\ge n_0. 
\end{align}

Now we consider the~representation 
\begin{equation}\label{6.19,1}
G_{\mu^{(n)}}(z)-G_{\nu^{(n)}}(z)=
\frac{r_{n1}(z)}{S_n(z)},\quad  z\in\Bbb C^+.
\end{equation}
The~relation (\ref{6.16,1}) leads to the~following estimate, for 
$z\in \Bbb C$ with $\Im z=1$,
\begin{align}\label{6.20,1}
\frac 12\big(1&+(|\Re z|-4)_+\big)\le |S(z)|
-\frac {6\cdot 10^3}{\sqrt {nm(z)}}\Big(|m_3(\mu)|+\frac {2m_4(\mu)}
{\sqrt n}\Big) \le|S_n(z)|\notag\\
&\le |S(z)|+\frac {6\cdot 10^3}{\sqrt {nm(z)}}\Big(|m_3(\mu)|+\frac {2m_4(\mu)}
{\sqrt n}\Big) \le 2\big(1+(|\Re z|-4)_+\big). 
\end{align}

Using (\ref{6.11b,1}), (\ref{6.19,1}), and (\ref{6.20,1}), we easily obtain
the~following inequality
\begin{align}\label{6.21,1}
\int\limits_{\Bbb R}|G_{\mu^{(n)}}(u+i)&-G_{\nu^{(n)}}(u+i)|\,du
\le\frac cn\Big(1+\frac{|m_3(\mu)|}{\sqrt n}+\frac{m_4(\mu)}n\Big)
\int\limits_{\Bbb R}\frac{du}{1+u^2}\notag\\
&\le\frac cn\Big(\frac{|m_3(\mu)|}
{\sqrt n}+\frac{m_4(\mu)}n\Big)\le \frac cn,\quad n\ge n_2, 
\end{align} 
and, for $ x\in[-2,2]$, using (\ref{6.11b,1}) and the~estimate
$|S_n(z)|\ge 1/10,\,z\in D_a$, we deduce
\begin{align}\label{6.22,1} 
\int\limits_a^1|G_{\mu^{(n)}}(x+iu)&-G_{\nu^{(n)}}(x+iu)|\,du\le
\frac c{n}\Big(1+\frac {|m_3(\mu)|}{\sqrt n}
+\frac {m_4(\mu)|\log a|}{n}\Big)\notag\\
&\le \frac c{n}\big(|m_3(\mu)|+(m_4(\mu))^{1/2}\big),\quad n\ge n_2. 
\end{align}

In order to prove the~upper bound of $\Delta(\mu^{(n)},\mu_w)$ for $n\ge n_2$
we apply Lemma~\ref{3.1a} with $v=a$. Since $m_4(\mu)<\infty$, 
it is well known that $m_2(\mu^{(n)})<\infty$ and
the~assumption (\ref{3.3a}) obviously holds.
Therefore Lemma~\ref{3.1a}, (\ref{6.17,1}), (\ref{6.18,1}), (\ref{6.21,1}),
and (\ref{6.22,1}) together imply the~estimate 
(\ref{2.9,a}). 

Hence, Theorem~\ref{2.3,1} is proved.
$\square$

{\bf Proof of Proposition~\ref{2.4,1}.}
Let $p,q>0$ and $p+q=1, p-q\ne 0$. We assume for definiteness that
$p-q<0$. Let $\nu$ be a~measure such that
$\nu(\{-p\})=q$ and $\nu(\{q\})=p$. It easy to see that $\nu_n=\mu_n$
for all $n=1,\dots$.
The~corresponding transforms 
are given by
$$
G_{\nu}(z)=\frac q{z+p}+\frac p{z-q} \quad\text{and}\quad
F_{\mu_n}(z)=F_{\nu_n}(z)=\frac{F_{\nu}(z\sqrt{npq})}{\sqrt{npq}}
=z-\frac 1n\cdot\frac 1{z+\tilde c/\sqrt n},
$$
where $\tilde c:=(p-q)/\sqrt{pq}$. A~simple calculations show that,
for $z\in\Bbb C^+$,
\begin{equation}\label{6.22,2}
\phi_{\mu_n}(z)=F_{\mu_n}^{(-1)}(z)-z=\frac 12\Big(-z-\frac {\tilde c}
{\sqrt n}+\sqrt{\Big(z-\frac {\tilde c}{\sqrt n}\Big)^2
+4\Big(\frac {\tilde c}{\sqrt n}z+\frac 1n\Big)}\Big),
\end{equation}
where we choose the~branch of the~square root which is positive for 
$z\ge 1/|\tilde c|$. Since $\phi_{\mu^{(n)}}(z)=n\phi_{\mu_n}(z)$,
we get $F_{\mu^{(n)}}^{(-1)}(z)=n\phi_{\mu_n}(z)+z$. Using this relation 
and (\ref{6.22,2}),we obtain with the~help of a~tedious but straightforward
calculation the~formula, for $z\in\Bbb C^+$, 
\begin{align}
F_{\mu^{(n)}}(z)&=\frac n{2(n-1)}\Big(\frac {n-2}nz-\frac{\tilde c}{\sqrt n}+
\sqrt{\Big(\frac {n-2}nz-\frac{\tilde c}{\sqrt n}\Big)^2-
4\frac{n-1}n\Big(1-z\Big(\frac{z}{n}+\frac{\tilde c}{\sqrt n}\Big)\Big)}
\,\Big)\notag\\
&=2\frac{1-z\Big(\frac{z}{n}+\frac{\tilde c}{\sqrt n}\Big)}
{\frac {n-2}nz-\frac{\tilde c}{\sqrt n}-\sqrt{\Big(\frac {n-2}nz
-\frac{\tilde c}{\sqrt n}\Big)^2-4\frac{n-1}n\Big(1-z\Big(\frac{z}{n}
+\frac{\tilde c}{\sqrt n}\Big)\Big)}}\notag
\end{align}
and hence, for the~same $z$,
\begin{align}
\Big(1-&z\Big(\frac{z}{n}+\frac{\tilde c}{\sqrt n}\Big)\Big) 
G_{\mu^{(n)}}(z)\notag\\
&=\frac 12\Big(\frac {n-2}nz-\frac{\tilde c}{\sqrt n}
-\sqrt{\Big(\frac {n-2}nz-\frac{\tilde c}{\sqrt n}\Big)^2-
4\frac{n-1}n\Big(1-z\Big(\frac{z}{n}+\frac{\tilde c}{\sqrt n}\Big)\Big)}
\,\Big)\notag\\
&=\frac 12\Big(\frac {n-2}nz-\frac{\tilde c}{\sqrt n}-
\sqrt{(z-x_1)(z-x_2)}\Big),\notag
\end{align}
where 
\begin{align}
x_1&:=-\frac{\tilde c}{\sqrt n}-
2\sqrt{1-\frac 1n}=-2-\frac{\tilde c}{\sqrt n}+\frac{2\theta}n,\notag\\
x_2&:=-\frac{\tilde c}{\sqrt n}+2\sqrt{1-\frac 1n}=2-\frac{\tilde c}
{\sqrt n}+\frac{2\theta}n.\notag
\end{align}
Here and below $\theta$ are real-valued quantities such that $|\theta|\le 1$.

Using the~Stieltjes-Perron inversion formula~(\ref{3.3}), we have
\begin{equation}\label{6.23,1}
\mu^{(n)}((-\infty,u))=\frac 1{2\pi}\int\limits_{x_1}^u p(x)\,dx,
\,\, x_1<u<x_2,\,\,\text{where}\,\, p(x):=
\frac{\sqrt{(x-x_1)(x_2-x)}}{1-x\Big(\frac{x}{n}+\frac{\tilde c}
{\sqrt n}\Big)}.
\end{equation}
Put $y=x+\tilde c/\sqrt n$. We note that
$$
(x-x_1)(x_2-x)=4-y^2-\frac 4n,\quad
\Big(1-x\Big(\frac{x}{n}+\frac{\tilde c}{\sqrt n}\Big)\Big)^{-1}=
1+\frac{\tilde c}{\sqrt n}y+\frac{\theta c(p)}n,
$$
for $x_1<x<x_2$. Here and below we shall denote by $c(p)$ positive constants 
depending on $p$ only. Since, for $-2+\frac {3|\tilde c|}{2\sqrt n}
\le x\le 2-\frac {3|\tilde c|}{2\sqrt n}$ and
for sufficiently large $n\ge n_0(p)$,
\begin{equation}
\sqrt{4-y^2-4/n}-\sqrt{4-y^2}=-\frac 4{\sqrt{4-y^2-4/n}+\sqrt{4-y^2}}\,\frac 1n=
\frac{c(p)\theta}{n^{3/4}},
\end{equation}
we easily obtain, for $x$ and $n$ as above,
$$
p(x)=\sqrt{4-y^2}+\frac{\tilde c y}{\sqrt n}\sqrt{4-y^2}
+\frac{c(p)\theta}{n^{3/4}}.
$$
Applying this formula to (\ref{6.23,1}) we deduce, 
for $-2\le u\le 0$,
\begin{align}
&\mu^{(n)}((-\infty,u))-\mu_w((-\infty,u))=\int\limits_u^{u+\tilde c/\sqrt n}
\frac{\sqrt{4-x^2}}{2\pi}\,dx-\frac{\tilde c(4-(u+\tilde c/\sqrt n)^2)^{3/2}}
{6\pi\sqrt n}+\frac{\theta c(p)}{n^{3/4}}\notag\\
&=\frac{\tilde c(4-(u+\theta\tilde c/\sqrt n)^2)^{1/2}}{2\pi\sqrt n}
-\frac{\tilde c(4-(u+\tilde c/\sqrt n)^2)^{3/2}}{6\pi\sqrt n}
+\frac{\theta c(p)}{n^{3/4}}.\notag
\end{align}
The~assertion of the~proposition now follows immediately from this relation.
$\square$.

\vspace {0,5 cm}

{\bf Proof of Theorem~\ref{2.5,1}.}
By Corollary~\ref{3.2a}, $G_{\mu^{(n)}}(z)=1/F_{\mu^{(n)}}(z),\,z\in \Bbb C^+$,
where $F_{\mu^{(n)}}(z):=F_{\mu_1}(Z_1(B_nz))/B_n=\dots
=F_{\mu_n}(Z_n(B_nz))/B_n$. In this formula
$Z_j(z),\,j=1,\dots,n$, are in the~class $\mathcal F$ and
are the~solutions of functional equations (\ref{3.7}). Without loss of 
generality, we assume that $\min_{j=1,\dots,n}m_2(\mu_j)\ge 1$ and
$\min_{j=1,\dots,n}m_2(\mu_j)=m_2(\mu_1)$. 
Denote $S_n(z):=Z_1(B_nz)/B_n$ and let, as in the~proof Theorem~\ref{2.3,1},
$S(z):=\frac 12(z+\sqrt{z^2-4})$. Note that $1/S_n(z)=G_{\nu^{(n)}}$
for some p-measure $\nu^{(n)}$. 

We prove the~inequality (\ref{2.10,b}) for $L_n\le c$ with a~sufficiently small
positive absolute constant $c$. For $L_n\ge c$ (\ref{2.10,b}) holds obviously.
From (\ref{3.7}) we have the~relation
\begin{equation}\label{6.4}
Z_1(z)-z=F_{\mu_2}(Z_2(z))-Z_2(z)+F_{\mu_3}(Z_3(z))-Z_3(z)+\dots
+F_{\mu_n}(Z_n(z))-Z_n(z),
\end{equation}
and
\begin{equation}\label{6.5}
F_{\mu_1}(Z_1(z))=F_{\mu_2}(Z_2(z))=\dots=F_{\mu_n}(Z_n(z)),
\quad z\in\Bbb C^+.
\end{equation}
By (\ref{6.4,1}), we note that
\begin{equation}\label{6.5a}
F_{\mu_j}(Z_j(z))-Z_j(z)=\frac {1-Z_j(z)G_{\mu_j}(Z_j(z))}
{Z_j(z)G_{\mu_j}(Z_j(z))}Z_j(z)=-\frac{r_{n,j}(z)}{1+r_{n,j}(z)}Z_j(z),
\quad z\in\Bbb C^+,
\end{equation}
where
\begin{equation}\label{6.5aa}
r_{n,j}(z):=\frac 1{Z_j(z)}\int\limits_{\Bbb R}\frac{u^2\,
\mu_j(du)}{Z_j(z)-u}=\frac{m_2(\mu_j)}{Z_j^2(z)}+\frac 1{Z_j^2(z)}
\int\limits_{\Bbb R}\frac{u^3\,\mu_j(du)}{Z_j(z)-u}.
\end{equation}
In addition, by (\ref{6.5}), we have
\begin{equation}\label{6.5b}
\frac{Z_1(z)}{Z_j(z)}=\frac{Z_1(z)G_{\mu_1}(Z_1(z))}
{Z_j(z)G_{\mu_j}(Z_j(z))}=\frac{1+r_{n,1}(z)}{1+r_{n,j}(z)},
\quad z\in\Bbb C^+.
\end{equation}

Since $\Im Z_j(B_nz)\ge B_n\Im z$, we obtain from 
(\ref{6.5aa}) that $|r_{n,j}(B_nz)|\le 1/10,\,j=1,\dots,n$, 
for $\Im z\ge c_3 M_n$, where $M_n:=(\max_{j=1,\dots,n} m_2(\mu_j))^{1/2}/
B_n$ and $c_3$ is a~sufficiently large absolute constant. 
Moreover, we deduce
from (\ref{6.5a}) and (\ref{6.5aa}) the~following estimates
\begin{align}\label{6.6a}
\Big|F_{\mu_j}(Z_j(B_nz))-Z_j(B_nz)&+\frac{m_2(\mu_j)}{Z_j(B_nz)}\Big|
\le\frac {\beta_3(\mu_j)}{|Z_j(B_nz)|B_n\Im z}
\notag\\
&+\frac {2m_2(\mu_j)}{|Z_j(B_nz)|^2B_n\Im z}\Big(m_2(\mu_j)
+\frac {\beta_3(\mu_j)}{B_n\Im z}\Big)
\end{align}
and
\begin{equation}\label{6.6}
\Big|F_{\mu_j}(Z_j(B_nz))-Z_j(B_nz)+\frac{m_2(\mu_j)}{Z_j(B_nz)}\Big|
\le\frac {\beta_3(\mu_j)}{|Z_j(B_nz)|B_n\Im z}+
\frac 2{|Z_j(B_nz)|^3}\Big(m_2(\mu_j)+\frac {\beta_3(\mu_j)}{B_n\Im z}\Big)^2
\end{equation}
for $\Im z\ge a_1:=c_3M_n$. 
In the~same way we obtain from (\ref{6.5b}) 
the~following inequalities
\begin{equation}\label{6.7}
\Big|\frac{Z_1(B_nz)}{Z_j(B_nz)}-1\Big|\le \frac{2}{B_n\Im z}
\Big(\frac {m_2(\mu_1)}{|Z_1(B_nz)|}+\frac{m_2(\mu_j)}{|Z_j(B_nz)|}\Big)
\le \frac 1{10}
\end{equation}
for $\Im z\ge a_1$ and $j=2,\dots,n$. Using (\ref{6.7}) we 
conclude that, for $\Im z\ge a_1$,
\begin{align}\label{6.8}
&\Big|\frac {m_2(\mu_2)}{Z_2(B_nz)}+\dots
+\frac {m_2(\mu_n)}{Z_n(B_nz)}-\frac{B_n^2-m_2(\mu_1)}{Z_1(B_nz)}
\Big|\notag\\
&\le \sum_{j=2}^n \frac{2m_2(\mu_j)}{|Z_1(B_nz)|B_n\Im z}\Big(\frac {m_2(\mu_1)}
{|Z_1(B_nz)|}+\frac{m_2(\mu_j)}{|Z_j(B_nz)|}\Big)
\le\frac{8}{|Z_1(B_nz)|^2B_n\Im z}\sum_{j=2}^n m_2^2(\mu_j).
\end{align}

In view of (\ref{6.6a}), (\ref{6.7}), and (\ref{6.8}),  
(\ref{6.4}) yields for $\Im z\ge a_1$ the~functional equation 
\begin{equation}\label{6.9}
S_n(z)-z=
-\frac {1-\widehat{r}_n(z)} {S_n(z)},
\end{equation}
where $\widehat{r}_n(z)$ is an~analytic function on $\Bbb C^+_{a_1}$
which admits the~upper bound 
\begin{equation}\notag
|\widehat{r}_n(z)|\le\frac 2{B_n^3\Im z}\sum_{j=1}^n\beta_3(\mu_j)+ 
\frac {12}{(B_n^2\Im z)^2}\sum_{j=1}^nm_2^2(\mu_j)
+\frac 4{B_n^5(\Im z)^3}\sum_{j=1}^nm_2(\mu_j)\beta_3(\mu_j)
+\frac {m_2(\mu_1)}{B_n^2}
\end{equation}
for $\Im z\ge a_1$. Using the~well-known inequalities
\begin{equation}\label{6.9a}
\frac 1{B_n^4}\sum_{j=1}^nm_2^2(\mu_j)\le\min\{M_n^2,L_n^{4/3}\},
\quad \frac 1{B_n^5}\sum_{j=1}^nm_2(\mu_j)\beta_3(\mu_j)\le M_n^2L_n,
\quad L_n\ge \frac 1{\sqrt n},
\end{equation}
we finally arrive at
\begin{equation}\label{6.10}
|\widehat{r}_n(z)|
\le\frac{2L_n}{\Im z}+\frac{12\min\{M_n^2,L_n^{4/3}\}}
{(\Im z)^2}+\frac{4M_n^2L_n}{(\Im z)^3}+L_n^2
\le \frac {20}{c_4}<\frac 1{10}
\end{equation}
for $\Im z\ge a_2:=c_4 (L_n+\min\{M_n,L_n^{2/3}\}+M_n^{2/3}L_n^{1/3})$, 
where $c_4>c_3$ is a~sufficiently large absolute constant. It follows
from (\ref{6.9}) and (\ref{6.10}) that 
\begin{equation}\label{6.10a}
10^{-1}\le |S_n(z)|\le 10,\qquad z\in D_{a_2},
\end{equation}
where the~closed domain $D_{a_2}$ is defined in the~proof of Theorem~2.4.
Using this inequality and (\ref{6.6})--(\ref{6.8}), we may improve
the~estimate (\ref{6.10}) for $z\in D_{a_2}$. Using as well (\ref{6.9a}) 
and the~well-known estimate 
$$
\frac 1{B_n^6}\sum_{j=1}^n\beta_3^2(\mu_j)\le L_n^2,
$$
we obtain the~following bound
\begin{align}\label{6.10b}
|\widehat{r}_n(z)|&\le\frac 2{B_n^3\Im z}\sum_{j=1}^n\beta_3(\mu_j)+ 
\frac {10^4}{B_n^4\Im z}\sum_{j=1}^nm_2^2(\mu_j)+\frac {10^4}
{B_n^6(\Im z)^2}
\sum_{j=1}^n\beta_3^2(\mu_j)+\frac {m_2(\mu_1)}{B_n^2}\notag\\
&\le 5\frac{L_n}{\Im z},
\quad z\in D_{a_2}.
\end{align}
By (\ref{6.10}), this estimate holds for $z\in\Bbb C^+$ such that $\Im z=1$.

Now we repeat the~arguments of the~proof of Theorem~\ref{2.3,1}.
Solving equation (\ref{6.9}) we see that
\begin{equation}\label{6.11}
S_n(z)=\frac 12\Big(z+\sqrt{\widehat{\rho}_n(z)}\Big),\quad
\Im z\ge a_2,
\end{equation} 
where $\widehat{\rho}_n(z):=z^2-4+4\widehat{r}_n(z)$. 

Write the~formula, for $z\in\Bbb C^+_{a_2}$,
\begin{equation}\label{6.12}
\frac 1{S_n(z)}-\frac 1{S(z)}=\frac{S(z)-S_n(z)}{S(z)S_n(z)}=
\frac 1{S(z)S_n(z)}\cdot\frac{\widehat{r}_n(z)}{\sqrt{z^2-4}+\sqrt{z^2-4
+4\widehat{r}_n(z)}}.
\end{equation}
Let $a_3:=3c_4L_n^{1/2}$. Note that the~well-known inequality
$M_n\le L_n^{1/3}$ implies $a_2<a_3$. 
Recalling that $|z^2-4|\ge m(z):=\max\{\Im z,((\Re z)^2-5)_+\},\,
0<\Im z\le 1$, we deduce from (\ref{6.10b}) that
$$
\Big|\frac{\widehat{r}_n(z)}{z^2-4}\Big|\le 
\frac{5L_n}{m(z)\Im z}\le\frac 1{10},\quad z\in D_{a_3}.
$$
Therefore we easily get, for $z\in D_{a_3}\cup \{z\in\Bbb C:\Im z=1\}$,
$$
|\sqrt{z^2-4}+\sqrt{z^2-4+4\widehat{r}_n(z)}|\ge \sqrt{|z^2-4|}\ge\sqrt{m(z)}.
$$  
Applying this estimate together with (\ref{6.10b}) to (\ref{6.12}), we  
conclude that, for $z\in D_{a_3}\cup\{z\in\Bbb C:\Im z=1\}$,
\begin{equation}\label{6.13}
\Big|\frac 1{S_n(z)}-\frac 1{S(z)}\Big|\le \frac{|\widehat{r}_n(z)|}
{|\sqrt{z^2-4}|}\frac 1{|S(z)||S_n(z)|}
\le 5\frac {L_n}{\sqrt{m(z)}\Im z|S(z)||S_n(z)|}.
\end{equation}


We conclude 
in the~same way as in (\ref{6.17,1}) and
(\ref{6.18,1}), using (\ref{6.13}), that is
\begin{equation}\label{6.14}
\int\limits_{\Bbb R}|G_{\mu_w}(u+i)-G_{\nu^{(n)}}(u+i)|\,du
\le cL_n\int\limits_{\Bbb R}\frac{du}{1+u^2}\le cL_n
\end{equation} 
and, for $x\in[-2,2]$,
\begin{equation}\label{6.15} 
\int\limits_{a_3}^1|G_{\mu_w}(x+iu)-G_{\nu^{(n)}}(x+iu)|\,du\le
c\int\limits_{a_3}^1
\frac {L_n}{u^{3/2}}\,du
\le c\frac {L_n}{a_3^{1/2}}\le cL_n^{3/4}. 
\end{equation}

Now we write
\begin{equation}\label{6.16}
G_{\mu^{(n)}}(z)-G_{\nu^{(n)}}(z)=
\frac{r_{n,1}(B_nz)}{S_n(z)},\quad  z\in\Bbb C^+.
\end{equation}
We deduce from (\ref{6.13}) the~following estimate, for 
$z\in D_{a_3}\cup\{z\in\Bbb C:\Im z=1\}$,
\begin{align}\label{6.17}
\frac 12\big(1&+(|\Re z|-4)_+\big)\le |S(z)|
-\frac {5L_n}{\sqrt {m(z)}\Im z}\le|S_n(z)|\notag\\
&\le |S(z)|+\frac {5L_n}{\sqrt {m(z)}\Im z}
\le 2\big(1+(|\Re z-4)_+\big). 
\end{align}
In addition we have, by (\ref{6.5aa}),
\begin{equation}\label{6.17a}
|r_{n,1}(B_nz)|\le \frac 1{(B_n|S_n(z)|)^2}\Big(m_2(\mu_1)+
\frac{\beta_3(\mu_1)}{B_n\Im z}\Big),\quad z\in D_{a_3}\cup\{z
\in\Bbb C:\Im z=1\}.
\end{equation}

Using (\ref{6.16})--(\ref{6.17a}), we easily obtain
the~following inequalities
\begin{equation}\label{6.18}
\int\limits_{\Bbb R}|G_{\mu^{(n)}}(u+i)-G_{\nu^{(n)}}(u+i)|\,du
\le c\Big(\frac 1n+L_n\Big)\int\limits_{\Bbb R}\frac{du}{1+u^2}
\le cL_n
\end{equation} 
and, for $x\in[-2,2]$,
\begin{equation}\label{6.19} 
\int\limits_{a_3}^1|G_{\mu^{(n)}}(x+iu)-G_{\nu^{(n)}}(x+iu)|\,du\le
c\Big(\frac 1n+L_n|\log a_3|\Big)\le cL_n |\log L_n|. 
\end{equation}

In order to prove the~upper estimate of $\Delta(\mu^{(n)},\mu_w)$
we apply again Lemma~\ref{3.1a} with $v=a_3$.
Lemma~{\ref{3.1a}, (\ref{6.14}), (\ref{6.15}), 
(\ref{6.18}), and (\ref{6.19}) together
imply the~estimate (\ref{2.10,b}) 
and the~theorem is proved.
$\square$
\vspace {0,5 cm}

{\bf Proof of Theorem~\ref{2.6,1}.}
For $k=1,\dots,n$ denote $\widehat{\mu}_{nk}((-\infty,x))
:=\mu_{k}((-\infty,nx+a_{nk})),\,x\in\Bbb R$, where 
$a_{nk}:=\int_{(-n,n)}u\mu_k(du)$. We shall now 
verify the~condition~(\ref{4.1}) with $k_n=n$ for the~measures 
$\widehat{\mu}_{nk}$. 
We obtain
\begin{align}
\varepsilon_{nk}=\int\limits_{\Bbb R}\frac {u^2}{1+u^2}\,
\widehat{\mu}_{nk}(du)
&=\int\limits_{\Bbb R}\frac {(u-a_{nk})^2}{n^2+(u-a_{nk})^2}\,
\mu_k(du)\notag\\
&\le\frac 1{n^2}\int\limits_{(-n,n)}(u-a_{nk})^2\mu_k(du)
+\int\limits_{\{|u|\ge n\}}
\mu_k(du).\notag
\end{align}
Therefore (\ref{4.1}) follows from (\ref{2.11,b}) and (\ref{2.11,c}).
Moreover, it follows from (\ref{2.11,b}) and (\ref{2.11,c}) that
\begin{equation}\label{6.6.1a}
\sum_{k=1}^n\varepsilon_{nk}\le \eta_{n1}+\eta_{n3}:=\eta_n.
\end{equation}

In the~proof of this theorem we use the~notation of Section~4 with
$k_n=n$ and $\tau=1$. 

From Corollary~3.4 we deduce the~relations (\ref{4.10}) and (\ref{4.11})
with $k_n=n$.
In addition $F_{\widehat{\mu}_n}(z)=F_{\widehat{\mu}_{n1}}(Z_{n1}(z)),\, 
z\in\Bbb C^+$, where $\widehat{\mu}_n:=\widehat{\mu}_{n1}\boxplus\dots\boxplus\widehat{\mu}_{nn}$.
By (\ref{4.9}) and (\ref{6.6.1a}), we get
\begin{equation}\label {6.6.1}
|\phi_{\widehat{\mu}_{n1}\boxplus\dots\boxplus\widehat{\mu}_{nn}}(z)|\le
|\phi_{\widehat{\mu}_{n1}}(z)|+\dots+|\phi_{\widehat{\mu}_{nn}}(z)|\le
c\sum_{k=1}^n\varepsilon_{nk}\le c\eta_n, 
\quad |z-i|\le 1/2.
\end{equation}

Since 
$$
\phi_{\widehat{\mu}_{n1}\boxplus\dots\boxplus\widehat{\mu}_{nn}}(z)
=(F_{\widehat{\mu}_{n1}}(Z_{n1}))^{(-1)}(z)-z=
Z_{n1}^{(-1)}(F_{\widehat{\mu}_{n1}}^{(-1)}(z))-z
$$ 
for $|z-i|\le 1/2$, we have, by (\ref{4.8}), the~relation
$$
\phi_{\widehat{\mu}_{n1}\boxplus\dots\boxplus\widehat{\mu}_{nn}}
(F_{\widehat{\mu}_{n1}}(z))=Z_{n1}^{(-1)}(z)-F_{\widehat{\mu}_{n1}}(z)
$$
for $|z-i|\le 1/4$.
Therefore we conclude by (\ref{4.8}) and (\ref{6.6.1})
that the~function 
$Z_{n1}^{(-1)}(z)$ is analytic in the~disk $|z-i|<1/4$ and 
$
|Z_{nk}^{(-1)}(z)-z)|\le c\eta_n
$
for $|z-i|< 1/4$. From this relation we see that 
\begin{equation}\label {6.6.2}
|Z_{n1}(z)-z|\le c\eta_n,\qquad |z-i|\le 1/8.
\end{equation} 
The~function $Z_{n1}(z)$ admits the~representation (\ref{4.14}). By
(\ref{6.6.2}), $|d_{n1}|\le c\eta_n$ and 
$\nu_{n1}(\Bbb R)\le c\eta_n$. Similar 
to (\ref{4.8}) we obtain
\begin{equation}\label{6.6.3}
|Z_{n1}(z)-z|\le c\eta_n\Big(1+\frac{1+|z|^2}{\Im z}\Big),
\quad z\in\Bbb C^+.
\end{equation}
Then we have, using (\ref{4.8}) and (\ref{6.6.3}), 
\begin{equation}\label{6.6.4}
|F_{\widehat{\mu}_{n1}}(Z_{n1}(z))-Z_{n1}(z)|\le  c\eta_n
\Big(1+\frac{1+|Z_{n1}(z)|^2}{\Im Z_{n1}(z)}\Big)
\le c\eta_n^{2/3}
\end{equation}
for $z=x+i\eta_{n}^{1/3},\,\eta_n^{1/6}\le x\le \eta_n^{1/6}$.
For such $z$ we finally get
\begin{equation}\label{6.6.4a}
|F_{\widehat{\mu}_n}(z)-z|\le c\eta_n^{2/3}.
\end{equation}
Since $F_{\widehat{\mu}_n}(z)\in\Cal F$ and therefore 
$|F_{\widehat{\mu}_n}(z)|\ge \Im z,\,z\in \Bbb C^+$, we conclude from
(\ref{6.6.4a}) that, for $z=x+i\eta_{n}^{1/3},\,
\eta_n^{1/6}\le x\le \eta_n^{1/6}$,
\begin{equation}\label{6.6.5}
\Big|G_{\widehat{\mu}_n}(z) -\ffrac 1z\Big|=\ffrac{|F_{\widehat{\mu}_n}(z)-z|}
{|F_{\widehat{\mu}_n}(z)||z|}\le c.
\end{equation}
From (\ref{6.6.5}) we get, for sufficiently large $n\ge n_4\ge c$,  
\begin{equation}\label{6.6.7}
-\frac 1{\pi}\int\limits_{\{|x|\le \eta_n^{1/6}\}}\Im G_{\widehat{\mu}_n}
(x+i\eta_n^{1/3})\,dx\ge\frac 1{\pi}\int\limits_{\{|x|\le \eta_n^{1/6}\}}
\frac{\eta_n^{1/3}}{x^2+\eta_n^{2/3}}\,dx-c\eta_n^{1/6}
\ge 1-c\eta_n^{1/6}. 
\end{equation}
On the~other hand we obtain
\begin{align}\label{6.6.8}
&-\frac 1{\pi}\int\limits_{|x|\le \eta_n^{1/6}}
\Im G_{\widehat{\mu}_n}(x+i\eta_n^{1/3})\,dx
=\frac 1{\pi}\int\limits_{\Bbb R}\Big(\arctan\frac{\eta_n^{1/6}-u}
{\eta_n^{1/3}}+\arctan\ffrac{\eta_n^{1/6}+u}{\eta_n^{1/3}}\Big)\,
\widehat{\mu}_n(du)\notag\\
&\le \widehat{\mu}_n(\{|u|\le 2\eta_n^{1/6}\})+1-\frac 1{\pi}
\arctan\eta_n^{-1/6}\le \widehat{\mu}_n(\{|u|\le 2\eta_n^{1/6}\})+c\eta_n^{1/6}.
\end{align}
From (\ref{6.6.7}) and (\ref{6.6.8}),  
for sufficiently large $n\ge n_4\ge c$, we have
$$
\widehat{\mu}_n(\{|u|\le 2\eta_n^{1/6}\})\ge 1-c\eta_n^{1/6}
$$
which immediately implies $L(\widehat{\mu}_n,\delta_0)\le c\eta_n^{1/6}$. 
By the~definition of $\mu^{(n)}$ and $\widehat{\mu}_n$,
we see that $L(\mu^{(n)},\widehat{\mu}_n)\le \eta_{n2}$. The~estimate 
(\ref{2.11,d}) is now an obvious consequence of the last two estimates.

Thus, Theorem~\ref{2.6,1} is proved. 
$\square$.

\end{document}